\documentclass[10pt,a4paper]{article}
\usepackage[latin1]{inputenc}
\usepackage{amsmath}
\usepackage{amsfonts}
\usepackage{amssymb}
\usepackage[english]{babel}
\usepackage{makeidx}
\usepackage{enumerate}
\usepackage{fancyhdr}
\usepackage[all]{xy}
\usepackage{lscape}
\usepackage{pict2e}

\usepackage{graphicx}
\graphicspath{{./}{SigPicture/}}
\DeclareGraphicsExtensions{.ps}

\newtheorem{theo}{Theorem}[section]
\newtheorem{prop}[theo]{Proposition}
\newtheorem{pte}[theo]{Property}
\newtheorem{defi}[theo]{Definition}
\newtheorem{lem}[theo]{Lemma}
\newtheorem{cor}[theo]{Corollary}

\newcommand{\eq}[1][r]
   {\ar@<-3pt>@{-}[#1]
    \ar@<-1pt>@{}[#1]|<{}="gauche"
    \ar@<+0pt>@{}[#1]|-{}="milieu"
    \ar@<+1pt>@{}[#1]|>{}="droite"
    \ar@/^2pt/@{-}"gauche";"milieu"
    \ar@/_2pt/@{-}"milieu";"droite"}
\newcommand{\citeg}[1]{\textbf{\cite{#1}}}
\newcommand{\refg}[1]{\textbf{\ref{#1}}}
\newcommand{\reffig}[1]{\textbf{Figure \ref{#1}}}

\newcommand{\RR}{\mathbb R}

\newcommand{\TT}{\mathcal T}
\newcommand{\Aa}{\mathcal A}
\newcommand{\EE}{\mathbb E}
\newcommand{\ZZ}{\mathbb Z}
\newcommand{\NN}{\mathbb N}

\newcommand{\BBB}{\mathcal{B}_0^c(\mathcal{T})}
\newcommand{\BB}{\mathcal{B}_0^c}
\newcommand{\B}{\mathcal{B}^c}
\newcommand{\ppp}{\mathfrak{p}_{0,\mathcal{T}}^c}

\newcommand{\p}{\mathfrak{p}}

\newcommand{\lambdo}{\lambda_0^c}
\newcommand{\inte}[1]{\overset{\!\!\!\! \circ}{#1}}

\title{\textbf{\huge{PV cohomology of pinwheel tilings, their integer group of coinvariants and gap-labelling}}}
\author{Ha\"ija MOUSTAFA}
\date{}

\begin{document}
\maketitle

\vspace{0cm}

\begin{abstract}
\noindent
In this paper, we first remind how we can see the "hull" of the pinwheel tiling as an inverse limit of simplicial complexes (\citeg{AndPut}) and we then adapt the PV cohomology introduced in \citeg{BelSav} to define it for pinwheel tilings. We then prove that this cohomology is isomorphic to the integer \v{C}ech cohomology of the quotient of the hull by $S^1$ which let us prove that the top integer \v{C}ech cohomology of the hull is in fact the integer group of coinvariants on some transversal of the hull. The gap-labelling for pinwheel tilings is then proved and we end this article by an explicit computation of this gap-labelling, showing that $\mu^t \big( C(\Xi,\ZZ) \big) = \dfrac{1}{264} \ZZ \left [ \dfrac{1}{5} \right ]$.
\end{abstract}

\vspace{0cm}

\tableofcontents

\newpage

\section{Introduction}

In this paper, we study some dynamical properties of the pinwheel tiling of the plane (the $(1,2)$-pinwheel tiling). Some remarks on possible extensions of this result will be made in conclusion.\\

\noindent
The study of tilings have gained in intensity since the discovery by physicists, in 1984, of a new material whose atomic distribution had forbidden symmetries for crystals (see \citeg{Sheetal}).\\
The atomic distribution wasn't the one of a crystal but it was very close to it : it was not periodic but nevertheless, it showed some order, it was "quasiperiodic". This material was called quasicrystals.\\
Quickly, mathematicians modeled such solids by aperiodic tilings and the physical properties of the material are closely related to the geometry of the tiling. This link was established by Jean Bellissard in \citeg{Bel1} and it is the content of the so called gap-labelling conjecture.\\
To state this conjecture, we need some definitions which will be given in more details in section \textbf{1.} of this paper.\\
First, to every tiling of the euclidean space, we can associate a topological space $\Omega$, called the continuous hull or the tiling space, which encode many properties of our tiling. This space is provided with an action of a subgroup $G$ of the isometries of the euclidean space (in the statement of the gap-labelling conjecture, $G$ will be the translations $\RR^n$) and thus we can consider the $C^*$-algebra associated to such dynamical system $(\Omega, G)$ which is the crossed product $C(\Omega) \rtimes G$.\\
Next, we assume that $\Omega$ is provided with an ergodic $G$-invariant probability measure $\mu$ which gives rise to a trace $\tau^\mu$ on $C(\Omega) \rtimes G$ and hence, to a linear map $\tau^\mu_* : K_0(C(\Omega) \rtimes G) \rightarrow \RR$ from the $K$-theory group of this $C^*$-algebra to the real numbers. The gap-labelling conjecture then predicts the image of $K_0(C(\Omega) \rtimes G)$ under this linear map.\\
Moreover, the hull contains a Cantor set $\Xi$, called the "canonical transversal" of $\Omega$, which is a sort of discretisation of the hull. The measure $\mu$ then induces a measure $\mu^t$ on $\Xi$ and the gap-labelling conjecture then expresses the link between the image of $K_0(C(\Omega) \rtimes \RR^n)$ under the trace and the image under $\mu^t$ of the integer valued functions on $\Xi$:\\

\noindent
\textbf{Conjecture : }  (\citeg{Bel1}, \citeg{BelHerZar}) \textit{
$$\tau_*^\mu \Big( K_0 \big( C(\Omega) \rtimes \RR^n \big) \Big) = \mu^t \big( C(\Xi,\ZZ) \big)$$
where $C(\Xi,\ZZ)$ is the space of continuous functions on $\Xi$ with values in $\ZZ$.}\\

\noindent
Since then, many works have been done to prove this conjecture.\\
First, the Pimsner-Voiculescu exact sequence gave the answer in dimension $1$ in \citeg{Bel1} and later, the conjecture was proved by van Elst in \citeg{vanElst} by iterating this exact sequence.\\
Using a spectral sequence, Bellissard, Kellendonk and Legrand proved the conjecture in dimension $3$ in \citeg{BelKelLeg}.\\
In $2002$, a general proof finally appears independently in several papers : by Bellissard, Benedetti and Gambaudo in \citeg{BelBenGam}, by Benameur and Oyono-Oyono in \citeg{BenOyo} and by Kaminker and Putnam in  \citeg{KamPut}.\\
The proof in \citeg{BenOyo} uses an index theorem for foliated spaces due to Alain Connes (see \citeg{Con1}) to link the analytical part $\tau_*^\mu \Big( K_0 \big(C(\Omega) \rtimes \RR^n \big) \Big)$ to a topological part $Ch_\tau \big( K_n \big( C(\Omega) \big) \big)$ which lies in $H^*_\tau(\Omega)$ the longitudinal cohomology group of $\Omega$ ($Ch_\tau$ is the longitudinal Chern character (see \citeg{MooSch})), this part being more computable.\\

\noindent 
By analogy, the gap labelling of the pinwheel tiling can be formulated as the computation of $\tau_*^\mu \Big( K_0 \big( C(\Omega) \rtimes \RR^2 \rtimes S^1 \big) \Big)$ in terms of the $\ZZ$-module of "patch frequencies" $\mu^t(C(\Xi,\ZZ))$.\\
We have adapted the method of Benameur and Oyono-Oyono in \citeg{Hai} for pinwheel tilings to obtain a similar result :\\

\noindent
\textbf{Theorem : } \textit{If $\Omega$ is the continuous hull of a pinwheel tiling, $\mu$ an ergodic invariant probability measure on $\Omega$ and $\mu^t$ the induced measure on the canonical transversal of $\Omega$, we have :
$$\tau_*^\mu \Big( K_0 \big( C(\Omega) \rtimes \RR^2 \rtimes S^1 \big) \Big) \subset \langle Ch_\tau \big( K_1( C(\Omega)) \big) , [C_{\mu^t}] \rangle$$
where $[C_{\mu^t}] \in H^\tau_3(\Omega)$ is the Ruelle-Sullivan current associated to the transverse measure $\mu^t$ (see \citeg{MooSch}) and $\langle \,, \, \rangle$ is the pairing of the longitudinal cohomology with the longitudinal homology.}\\

\noindent
Since the longitudinal Chern character factorizes through the usual Chern character in \v{C}ech cohomology (see section \textbf{3.}) and since the Ruelle-Sullivan current only sees the $H^3_\tau(\Omega)$ part, we will study, in this paper, the top integer \v{C}ech cohomology group of $\Omega$.\\
The aim of this paper is to study carefully the image under the Ruelle-Sullivan current of the top dimensional longitudinal Chern character and to relate it to the module of patch frequencies in order to solve completely the gap-labelling conjecture for the pinwheel tiling.\\ 
\\
\noindent
The structure of this paper is then the following: in section \textbf{2.} we remind some classical definitions in tiling theory. In particular, we remind the construction of the pinwheel tiling given by Radin in \citeg{Rad} and then we introduce the notion of continuous hull of pinwheel tiling, enumerating its properties. We next turn to the definition of the canonical transversal of the hull which allows us to see the continuous hull of pinwheel tiling as a foliated space in a well known way.\\

\noindent
In section \textbf{3.}, a trick shows that $\check{H}^3(\Omega ; \ZZ)$ is isomorphic to $\check{H}^2(\Omega/S^1 ; \ZZ)$.\\
We then study the top integer \v{C}ech cohomology of $\Omega/S^1$.\\
To do this, we are using results and ideas developed in several papers (\citeg{AndPut}, \citeg{BelSav}).\\
First, we use the idea, initiated by Anderson and Putnam in their paper \citeg{AndPut}, to see $\Omega/S^1$ as an inverse limit of homeomorphic simplicial complexes.\\
Specifically, as the pinwheel tiling isn't "forcing its border", we will use the collared version of their construction.\\
This allows to use simplicial methods to compute $\check{H}^2(\Omega/S^1 ; \ZZ)$.\\
This was extended by Bellissard and Savinien in \citeg{BelSav} to compute the cohomology of tilings in term of the PV cohomology of its prototile space.\\
In this section, we adapt their method to prove 
$$\check{H}^* \big( \Omega/S^1 ; \ZZ \big) \simeq H^*_{PV} \big( \B_0 ; C(\Xi_\Delta^2,\ZZ) \big).$$
\noindent
where $H^*_{PV}$ is the PV cohomology and $\Xi_\Delta^2$ is a new transversal .\\
\\
\noindent
The interesting point in this new cohomology is that the cochains are in fact classes of continuous functions with integer values on the transversal $\Xi_\Delta^2$ of the hull which is a first step toward the module of "patch frequencies" of the pinwheel tiling, related to the continuous functions with integer values on the canonical transversal.\\

\noindent
The key point established in this section is that the top PV cohomology of the pinwheel tiling is isomorphic to the integer group of coinvariants of the transversal $\Xi_\Delta^2$ (this notion of integer group of coinvariants is given in this section, it's the quotient of the continuous function on $\Xi_\Delta^2$ by the "local" coinvariants, in a similar way to the definition given in \citeg{Kelcoinv}) : \\

\noindent
\textbf{Theorem 3.16 :} \textit{
$$H^2_{PV}(\B_0;C(\Xi_\Delta,\ZZ)) \cong C(\Xi^2_\Delta,\ZZ)/ H_{\Xi_\Delta^2}$$
where $H_{\Xi_\Delta^2}$ is a subgroup of $C(\Xi^2_\Delta,\ZZ)$ such that for all $h \in H_{\Xi_\Delta^2}$, $\mu_2^t(h)=0$, for the measure $\mu^t_2$ induced by $\mu$ on $\Xi_\Delta^2$.}\\

\noindent
which leads to the important corollary that the top integer \v{C}ech cohomology group of the hull is isomorphic to the integer group of coinvariants :\\

\noindent
\textbf{Corrolary 3.17 :} \textit{The top integer \v{C}ech cohomology of the hull is isomorphic to the integer group of coinvariants of $\Xi_\Delta^2$ :
$$\check{H}^3 \big( \Omega ; \ZZ \big) \simeq C(\Xi_\Delta^2,\ZZ) / H_{\Xi_\Delta^2}.$$}

\noindent
In section \textbf{4.}, this theorem associated to the study of the image under the Ruelle-Sullivan map of the top cohomology of $\Omega$ gives the desired  gap-labelling of the pinwheel tiling :\\

\noindent
\textbf{Theorem 4.1 :} \textit{
If $\TT$ is a pinwheel tiling, $\Omega = \Omega(\TT)$ its hull provided with an invariant ergodic probability measure $\mu$ and $\Xi$ its canonical transversal provided with the induced measure $\mu^t$, we have :
$$\tau_*^\mu \Big( K_0 \big( C(\Omega) \rtimes \RR^2 \rtimes S^1 \big ) \Big) = \mu^t \big( C(\Xi,\ZZ) \big).$$
}

\noindent
We finally end section \textbf{4.} by an explicit computation of this image.\\
Viewing $C(\Xi,\ZZ)$ as a direct limit and exhausting the collared prototiles of the pinwheel tiling, we then prove that, thanks to a result in \citeg{Eff}, we have :

$$\tau_*^\mu \Big( K_0 \big( C(\Omega) \rtimes \RR^2 \rtimes S^1 \big) \Big) = \frac{1}{264} \ZZ \left [ \frac{1}{5} \right].$$

\noindent
This result shows that the gap-labelling of the pinwheel tiling is given by the $\ZZ$-module of its "patch frequencies".\\
\\
A natural question is whether this is a general fact.\\
\\
\noindent
\textbf{Aknowledgements. } It is a pleasure for me to thank my advisor Herv\'e Oyono-Oyono who always supported and advised me during this work.\\
I also want to thank Jean Bellissard for useful discussions on the gap-labelling conjecture, Ian Putnam and Michael Whittaker for useful discussions on pinwheel tilings and their gap-labelling.\\
I am also grateful to Dirk Frettl\"oh for useful conversations on diffraction and patch frequencies.\\
I am also indebted to the SSM department of Victoria BC University for its hospitality during my stay in march $2009$ where this paper was written.\\

\newpage

\section{Reminders}

\subsection{Pinwheel tiling and continuous hull}

\vspace{0.5cm}
A \textbf{tiling of the plane} is a countable family $P=\{t_1, t_2 \ldots \}$ of non empty compact subsets $t_i$ of $\RR^2$, called \textbf{tiles} (each tile being homeomorphic to the unit ball), such that:
\begin{itemize}
\item[$\bullet$] ${\displaystyle \bigcup_{i \in \mathbb{N}} t_i = E_2}$ where $E_2$ is the euclidean plane with a fixed origin $O$;
\item[$\bullet$] Tiles meet each other only on their border ;
\item[$\bullet$] Tiles's interiors are pairwise disjoint.\\
\end{itemize} 

\noindent
We are interested in the special case where there exists a finite family of tiles $\{p_1, \ldots, p_n \}$, called \textbf{prototiles}, such that each tile $t_i$ is the image of one of these prototiles under a \textit{rigid motion} (i.e. a direct isometry of the plane).\\
In fact this paper will focus on the particular tiling called \textbf{pinwheel tiling} or \textbf{(1,2)-pinwheel tiling} which is obtained by a substitution explained below.\\

\begin{figure}[ht]
\begin{center}
\includegraphics[width=8cm]{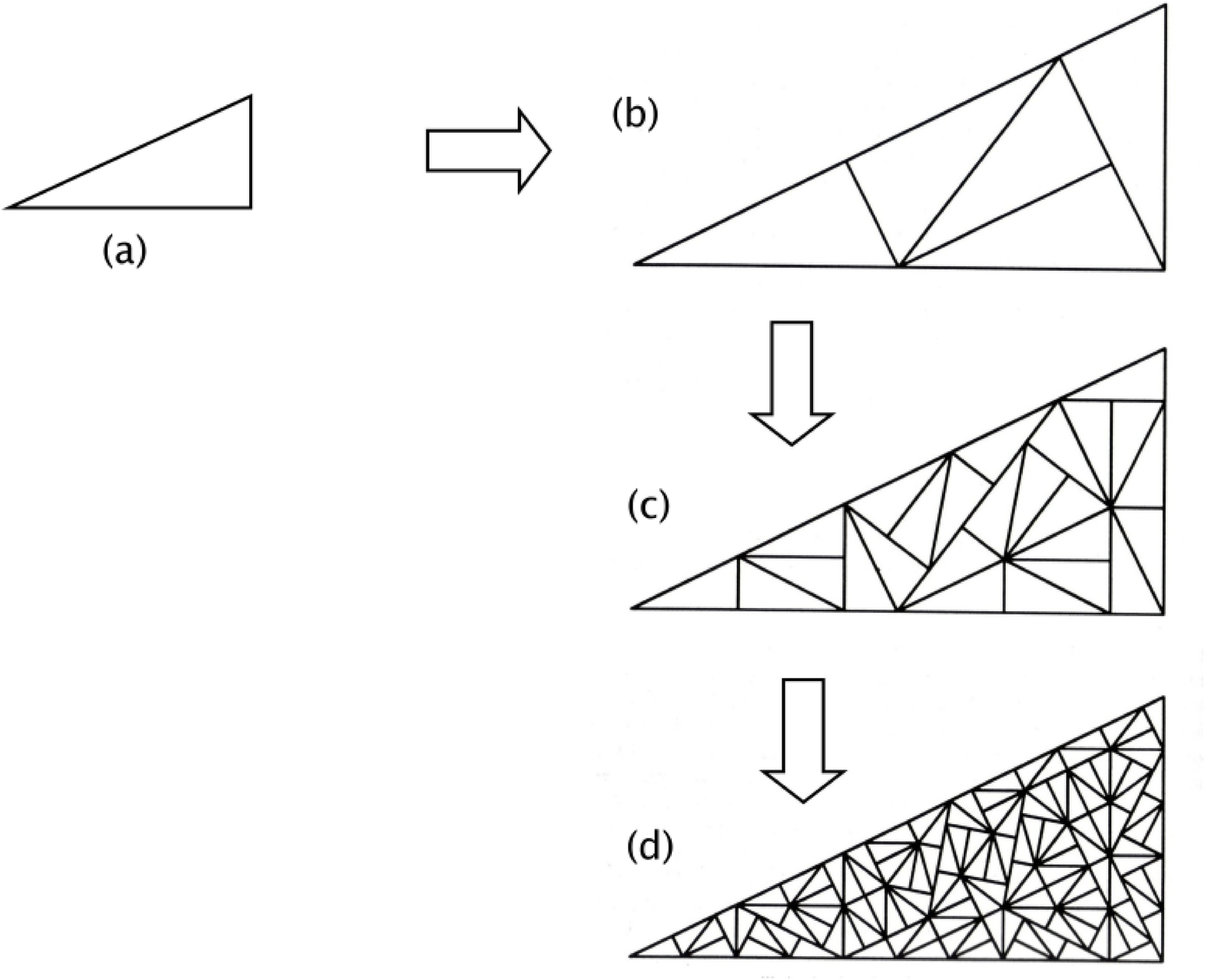}
\end{center}
\caption{Substitution of the pinwheel tiling.}
\label{Pin}
\end{figure}

\noindent
Our construction of a pinwheel tiling is based on the construction made by Charles Radin in \citeg{Rad}. It's a tiling of the plane obtained by the substitution described in \reffig{Pin}.\\

\noindent
This tiling is constructed from two prototiles, the right triangle in \reffig{Pin}\textbf{.(a)} with legs $1$, $2$ and $\sqrt{5}$ and its mirror image.\\

\begin{figure}[ht]
\begin{center}
\includegraphics[width=\textwidth]{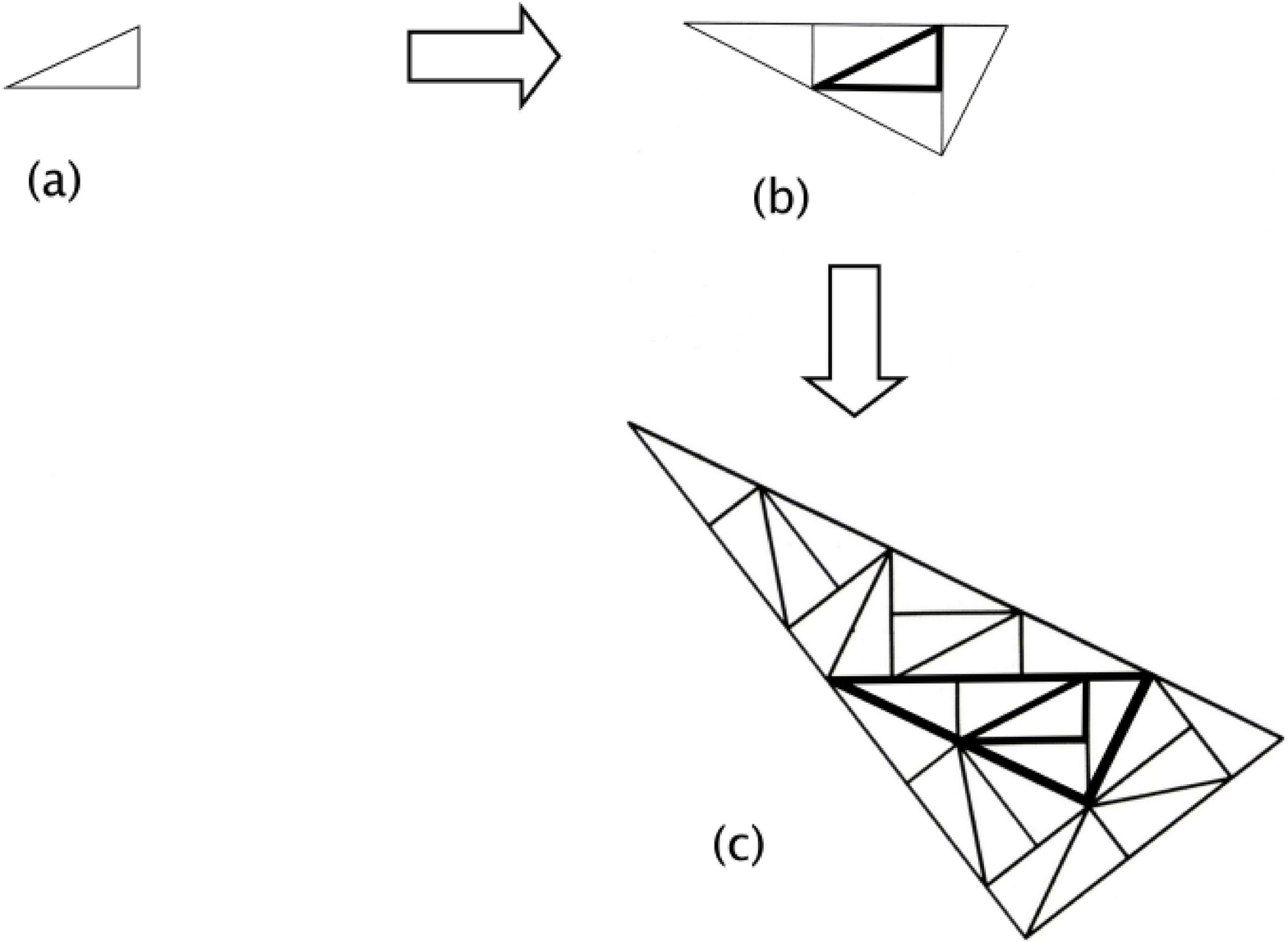}
\end{center}
\caption{Construction of a pinwheel tiling.}
\label{pinwrot}
\end{figure}

\noindent
To obtain this tiling, we begin from the right triangle with the following vertices in the plane : $(0,0)$ , $(2,0)$ and $(2,1)$.\\
This tile and its reflection are called \textbf{supertiles of level 0} or \textbf{0-supertiles}.\\
We will next define \textbf{1-supertiles} as follows : take the right triangle with vertices $(-2,1)$, $(2,-1)$ and $(3,1)$ and take the decomposition of \reffig{Pin}\textbf{.(b)}. This $1$-supertile is thus decomposed in five $0$-supertiles, which are isometric copies of the first tile, with the beginning tile in its center (see \reffig{pinwrot}\textbf{.(b)}).

\noindent
We next repeat this process by mapping this $1$-supertile in a $2$-supertile with vertices $(-5,5)$, $(1,-3)$ and $(5,0)$ (see \reffig{pinwrot}\textbf{.(c)}).\\
Including this $2$-supertile in a $3$-supertile with correct orientation and so on, this process leads to the desired \textbf{pinwheel tiling} $\TT$.\\

\noindent
We will now attach to this tiling a topological space reflecting the combinatorial properties of the tiling into topological and dynamical properties of this space.\\

\noindent
For this, we observe that the direct isometries of the plane are acting on the euclidean plane $E_2$ where we have fixed the origin $O$.\\
Direct isometries $\EE^2=\RR^2 \rtimes SO(2)$ thus act naturally on our tiling $\TT$ on the right .\\
If $R_\theta$ denotes the rotation about the origin with angle $\theta$ and $s \in \RR^2$, $\TT.(s,R_\theta):=R_{-\theta}(\TT-s)$. 
We will also denote $(s,R_\theta)$ by $(s,\theta)$.

\begin{defi} \text{ }\\
A \textbf{patch} is a finite union of tiles of a tiling.\\
A tiling $\TT'$ is of \textbf{finite $\mathbf{\EE^2}$-type} or of \textbf{Finite Local Complexity} (\textbf{FLC}) if for any $R > 0$, there is only a finite number of patches in $\TT'$ of diameter less than $R$ up to direct isometries.\\
A tiling $\TT'$ of finite $\mathbf{\EE^2}$-type is \textbf{$\mathbf{\EE^2}$-repetitive} if, for any patch $\Aa$ in $\TT'$, there is $R(\Aa)>0$ such that any ball of radius $R(\Aa)$ intersects $\TT'$ in a patch containing a $\EE^2$-copy of $\Aa$.
\end{defi}

\noindent
The tiling $\TT$ is of finite $\EE^2$-type, $\mathbf{\EE^2}$-repetitive and non periodic for translations (see \citeg{Pet}).\\

\noindent
To attach a topological space to $\TT$, we define a metric on $\TT.\EE^2$ :\\
If $\TT_1$ and $\TT_2$ are two tilings in $\TT.\EE^2$, we define 
$$\hspace*{-2cm} A=\left \{\varepsilon \in \Big [ 0,\tiny{\frac{1}{\sqrt{2}}} \Big ] \; / \; \exists s,s' \in B^2_\varepsilon(0) \, , \, \theta,\theta' \in B^1_\varepsilon(0) \normalsize{\text{ s.t. }} \right.$$ 
$$\hspace*{3cm} \left. \TT_1.(s,\theta) \cap B_{\frac{1}{\varepsilon}}(O) =  \TT_2.(s',\theta') \cap B_{\frac{1}{\varepsilon}}(O) \right \}$$
where $B_{\frac{1}{\epsilon}}(O)$ is the euclidean ball centered in $O$ with radius $\frac{1}{\epsilon}$ and $B_\epsilon^i(0)$ are the euclidean balls in $\RR^i$ centered in $0$ and with radius $\epsilon$ (i.e. we consider direct isometries near $Id$).\\
Then, define :
$$d(\TT_1,\TT_2) = \left\{
 \begin{array}{cc}
Inf A & \text{ if } A \neq \emptyset\\
\frac{1}{\sqrt{2}} & \text{ else } \quad
\end{array} \right..$$
$d$ is a bounded metric on $\TT.\EE^2$. For this topology, a base of neighborhoods is defined by: two tilings $\TT_1$ and $\TT_2$ are close if, up to a small direct isometry, they coincide on a large ball around the origin.\\

\begin{defi} \text{ }\\
The \textbf{continuous hull} of $\TT$ is then the completion of $(\TT.\EE^2,d)$ and will be denoted $\Omega(\TT)$.
\end{defi}

\noindent
Let's enumerate some well known properties of this continuous hull:

\begin{pte}  (\citeg{KelPut}, \citeg{BelBenGam}, \citeg{BenGam} et \citeg{Rad}) \label{decompo}
\begin{itemize}
	\item $\Omega(\TT)$ is formed by finite $\EE^2$-type, $\mathbf{\EE^2}$-repetitive and non periodic (for translations) tilings and each tiling of $\Omega(\TT)$ has the same patches as $\TT$.
	\item $\Omega(\TT)$ is a compact space since $\TT$ is of finite $\EE^2$-type.
	\item Each tiling in $\Omega(\TT)$ are uniquely tiled by $n$-supertiles, for all $n \in \NN$.
	\item The dynamical system $(\Omega(\TT), \EE^2)$ is minimal since $\TT$ is repetitive, i.e each orbit under direct isometries is dense in $\Omega(\TT)$.
\end{itemize}
\end{pte}

\noindent
The last property of $\Omega(\TT)$ allows us to write $\Omega$ without mentioning the tiling $\TT$ (in fact, if $\TT' \in \Omega(\TT)$, $\Omega(\TT') = \Omega(\TT)$).

\begin{defi}
Any tiling in $\Omega$ is called a \textbf{pinwheel tiling}.
\end{defi}

\noindent
\textbf{Remark :} we can easily see that our continuous hull is the compact space $X_\phi$ defined by Radin and Sadun in \citeg{RadSad}.

\newpage

\subsection{The canonical transversal}

\vspace{0.5cm}
In this section, we will construct a Cantor transversal for the action of $\EE^2$ and we show that this transversal gives the local structure of a foliated space.\\
For this, we fix a point in the interior of the two prototiles of the pinwheel tiling. This, in fact, gives for any tiling $\TT_1$ in $\Omega$ (i.e constructed by these two prototiles), a punctuation of the plane denoted $\TT_1^{punct}$.\\
Define then $\Omega_0$ to be the set of every tilings $\TT_1$ of $\Omega$ such that $O \in \TT_1^{punct}$.\\
The \textbf{canonical transversal} is the space $\Omega_0 / SO(2)$.\\

\noindent
We can identify this space with a subspace of $\Omega$ by constructing a continuous section $s:\Omega_0/SO(2) \longrightarrow \Omega$.\\
To obtain such a section, we fix an orientation of the two prototiles of our tilings once for all. Hence when we consider a patch of a tiling in the transversal $\Omega_0$, there is only one orientation of this patch where the tile containing the origin have the orientation chosen for the prototiles.\\
Let then $[\omega] \in \Omega_0 / SO(2)$, there is only one $\theta \in [0;2\pi[$ such that the tile in $R_\theta(\omega)$ containing the origin has the good orientation.\\
We define $s([\omega]) := R_\theta(\omega)$.\\
$s$ is well defined because $\theta$ depends on the representative $\omega$ chosen but not $R_\theta(\omega)$.\\
$s:\Omega_0 / SO(2) \longrightarrow s(\Omega_0 / SO(2))$ is then a bijection. We easily see that $s$ is continuous and thus it is a homeomorphism from the canonical transversal onto a compact subspace $\Xi$ of $\Omega$.\\
We also call this space the \textbf{canonical transversal}.\\
We can see $\Xi$ as the set of all the tilings $\TT_1$ in $\Omega$ with the origin on the punctuation of $\TT_1$ and with the tile containing the origin in the orientation chosen for the prototiles.\\

\noindent
We then have :

\begin{prop} (\citeg{BenGam})\\
The canonical transversal is a Cantor space.
\end{prop}
 
\noindent
A base of neighborhoods is obtained as follows : consider $\TT' \in \Xi$ and $\Aa$ a patch around the origin in $\TT'$ then 
$$U(\TT',\mathcal{A}) = \{\TT_1 \in \Xi \; \vert \; \TT_1 = \TT'  \text{ on } \Aa \}$$
is a closed and open set in $\Xi$, called a \textbf{clopen} set.\\

\noindent
Before defining the foliated stucture on $\Omega$, we must study the rotations which can fix tilings in $\Omega$.\\
In pinwheel tilings, we can sometimes find regions tiled by supertiles of any level and so we introduce the following definition:

\begin{defi}
A region of a tiling which is tiled by $n$-supertiles for all $n \in \NN$ is called an \textbf{infinite supertile} or \textbf{supertile of infinite level}.
\end{defi}

\noindent
If a ball in a tiling $\TT_1$ fails to lie in any supertile of any level n, then $\TT_1$ is tiled by two or more supertiles of infinite level, with the offending ball straddling a boundary.\\
We can, in fact, construct a pinwheel tiling with two half-planes as infinite supertiles as follows: Consider the rectangle consisting of two $(n-1)$-supertiles  in the middle of a $n$-supertile. For each $n \geqslant 1$, orient this rectangle with its center at the origin and its diagonal on the $x$-axis, and fill out the rest of a (non-pinwheel) tiling $\TT_n$ by periodic extension. By compactness this sequence has a convergent subsequence, which will be a pinwheel tiling and which will consist of two infinite supertiles (this example comes from \citeg{RadSad}).\\
\\
Note that the boundary of an infinite supertile must be either a line, or have a single vertex, since it is tiled by supertiles of all levels.\\
We call such a line a \textbf{fault line}.

\begin{lem}
If $(s,\theta)$ fixes a pinwheel tiling $\TT'$ then $\theta \in \{ 0,\pi \} mod(2\pi)$.\\
Moreover, if $\theta = 0$ then $s=0$. In other terms, translations can't fix a pinwheel tiling.
\end{lem}

\noindent
\textbf{Proof :} Let's consider the different cases:
\begin{enumerate}
\item First, if the tiling $\TT'$ which is fixed by $(s,\theta)$ have no fault line (i.e have no infinite supertile), then $s=0$ and $\theta=0 (mod 2 \pi)$.\\
Indeed, let $x \in E_2$ be such that $\overrightarrow{Ox} = s$ then $O$ and $x$ is in the interior of a $m$-supertile since there isn't infinite supertiles (see p.29 in \citeg{RadSad}).\\
As no direct isometry fixes our prototiles, $s$ and $\theta$ must be zero.

\item Let's see the case in which $\TT'$ have some infinite supertiles.\\
By \citeg{RadSad} p.30, the number of infinite supertiles in $\TT'$ is bounded by a constant $K$ (in fact for pinwheel tilings, we can take $K=\frac{2\pi}{\alpha}$ where $\alpha$ is the smallest angle in the prototiles).\\
Thus, $\TT'$ doesn't contain more than $K$ infinite supertiles and in fact, it has only a finite number of fault lines. This will give us the result.\\
Indeed, since $(s,\theta)$ fixes $\TT'$, if $F$ is a fault line, $(s,\theta)$ sends it on another fault line $F_1$ in $\TT'$ and thus, $R_\theta$ sends $F$ on a line  parallel to $F_1$.\\
As there is only a finite number of fault lines in $\TT'$, there is $M \in \NN^*$, $m \in \ZZ^*$ such that $M \theta = 2 \pi m$.\\
If we now use results obtained in \citeg{RadSad} p.32, $\theta$ must be in the group of relative orientations $G_{RO}(Pin)$ of pinwheel tilings which is the subgroup of $SO(2)$ generated by $\frac{\pi}{2}$ and $2\alpha$.\\
Hence, if $\theta=2k\alpha + l \frac{\pi}{2}$ with $k \in \ZZ^*$ and $l \in \ZZ$, it would mean that $\alpha$ is rationnal with respect to $\pi$, which is impossible (see \citeg{Rad} p.664) hence $k=0$ and $\theta \in \left \{ 0, \, \dfrac{\pi}{2}, \, \pi, \, \dfrac{3\pi}{2} \right \}$ $mod(2 \pi)$.\\
Now, if we study the first vertex coronas (i.e the minimal patches around a vertex or around the middle point of the hypothenuses), there is only patches with a $2$-fold symmetry (\citeg{SadPriv} or see \reffig{prototuilescouronnees} and \reffig{prototuilescouronnees2} p.\pageref{prototuilescouronnees} and p.\pageref{prototuilescouronnees2}).\\
Thus, $\theta \in \left \{ 0, \, \pi \right \}$ $mod(2 \pi)$ and if $\theta = 0$, $s=0$ since pinwheel tilings are not fixed by translations.
\end{enumerate}

\begin{flushright}
$\square$
\end{flushright}

\noindent
We note that, in fact, there exists only $6$ pinwheel tilings with a $2$-fold symmetry up to rotations (\citeg{SadPriv}).\\
Hence, there is only $6$ orbits with fixed points for the $\RR^2 \rtimes SO(2)$ action on $\Omega$.\\
Moreover, there is exactly $6$ circles $F_1, \ldots, F_6$ containing fixed points for the $SO(2)$ action on $\Omega$ (of course, therefore, the $6$ orbits of these circles contain all the fixed points of the $\RR^2 \rtimes SO(2)$-action).\\

\noindent
We thus obtain the following important result on the dynamic of our tiling space:

\begin{theo} (\citeg{BenGam})
The continuous hull is a minimal foliated space.
\end{theo}

\noindent
\textbf{Proof :}
\begin{itemize}
\item[] The proof follows the one in \citeg{BenGam} except that, locally, $\Omega$ looks like an open subset of $SO(2)$ $\times$ an open subset of $\RR^2$ $\times$ a Cantor set instead of $SO(2)$ $\times$ an open subset of $\RR^2$ $\times$ a Cantor set, like in \citeg{BenGam}.\\

\noindent
$\Omega$ is covered by a finite number of open sets $U_i = \phi_i (V_i \times T_i)$ where :
\begin{itemize}
\item[$\bullet$] $T_i$ is a clopen set in $\Xi$;
\item[$\bullet$] $V_i$ is an open subset of $\RR^2 \rtimes SO(2)$ which read $V_i = \Gamma_i \times W_i$ where $W_i$ is an open subset of $\RR^2$ and $\Gamma_i$ an open subset of $SO(2)$ of the form $]l\pi/2-\pi/3 ; l\pi/2+\pi/3[$, $l \in \{0,1,2,3\}$;
\item[$\bullet$] $\phi_i:V_i \times T_i \longrightarrow \Omega$ is defined by $\phi_i(v,\omega_0) = \omega_0.v$.
\end{itemize}

\noindent
As we can find finite partitions of $\Xi$ in clopen sets with arbitrarily small diameter, it is possible to choose this
diameter small enough so that:

\begin{itemize}
\item[$\bullet$] the maps $\phi_i$ are homeomorphisms on their images;
\item[$\bullet$] whenever $\TT_1 \in U_i \cap U_j$, $\TT_1 = \phi_i(v,\omega_0) = \phi_j(v',\omega_0')$, the element $v'.v^{-1}$ is independent  of the choice of $\TT_1$ in $U_i \cap U_j$, we denote it by $g_{ij}$.
\end{itemize}
The transition maps read : $(v',\omega_0') = (g_{ij}.v,\omega_0.g_{ij}^{-1})$.\\
It follows that the boxes $U_i$ and charts $h_i = \phi_i^{-1} : U_i \longrightarrow V_i \times T_i$ define a foliated structure on $\Omega$.\\
By construction, the leaves of $\Omega$ are the orbits of $\Omega$ under the $\EE^2$-action.
\end{itemize}

\begin{flushright}
$\square$
\end{flushright}

\noindent
We must do several remarks now on the actions.\\
$\EE^2$ isn't acting freely on $\Omega$, even if the translations are, but we could adapt results of Benedetti and Gambaudo obtained in their paper \citeg{BenGam} studying the possible symmetries in our pinwheel tilings.\\
The $\EE^2$-action is not free on $\Omega_0$ too but the $SO(2)$-action is.\\

\noindent
Using the group of relative orientations $G_{RO}(Pin)$, we can see that each $\RR^2$-orbit of $\Omega$ is in fact a dense subset of $\Omega$ (see \citeg{HolRadSad}).

\newpage

\section{PV cohomology of pinwheel tilings and the integer group of coinvariants}

In \citeg{Hai}, we have obtained 
$$\tau^*_\mu\left( K_0 \big(C(\Omega) \rtimes \RR^2 \rtimes SO(2) \big) \right) \subset [C_{\mu^t}] \Big( Ch_\tau (K_1(C(\Omega)) \Big).$$

\noindent
There is a natural map (\citeg{MooSch}) $r^*:\check{H}^3(\Omega;\RR) \longrightarrow H^3_\tau(\Omega)$ obtained by inclusion of the sheaf $\mathcal{R}$ of germs of locally constant real-valued functions into the sheaf $\mathcal{R}_\tau$ of germs of continuous real-valued tangentially locally constant functions and since $\ZZ \subset \RR$, we also have a natural map, also denoted $r^*$, $r^*:\check{H}^3(\Omega;\ZZ) \longrightarrow H^3_\tau(\Omega)$.\\
We then have a factorization : $Ch_\tau = r^* \circ Ch$ where $Ch$ is the Chern character $Ch:K_1(\Omega) \longrightarrow \check{H}^{odd}(\Omega,\ZZ)$ (the Chern character can be defined with value in the integer odd \v{C}ech cohomology because we are in dimension 3).\\
The Ruelle-Sullivan current $[C_{\mu^t}]$ only takes into account the $H^3_\tau(\Omega)$ part of $Ch_\tau$ and thus we must focus on the top \v{C}ech cohomology $\check{H}^3(\Omega;\ZZ)$.\\
In fact, we will study $\check{H}^2(\Omega/S^1 ; \ZZ)$ since 
$$\check{H}^3(\Omega; \ZZ) \simeq \check{H}_c^3(\Omega \setminus F; \ZZ) \simeq \check{H}_c^2 \big( (\Omega \setminus F)/S^1 ; \ZZ) \simeq \check{H}^2(\Omega/S^1; \ZZ).$$
The left hand side and the right hand side isomorphisms are obtained by the long exact sequence in cohomology relative to the pairs $(\Omega,F)$ and  $(\Omega/S^1,F/S^1)$ and use the fact that $F$ is of dimension 1 and $F/S^1$ of dimension 0.\\
The isomorphism in the middle is obtained by a Gysin sequence since the projection $\Omega \setminus F \longrightarrow (\Omega \setminus F)/S^1$ is a $S^1$-principal bundle as the $S^1$-action is free on $\Omega \setminus F$ (see \citeg{Bred}).\\

\noindent
To study this cohomology, we use techniques developed in the paper of Bellissard and Savinien (\citeg{BelSav}) to show that the top integer \v{C}ech cohomology of $\Omega / S^1$ is isomorphic to the integer group of coinvariants associated with a certain transversal.\\

\noindent
This will be achieved by using an idea first introduced by Anderson and Putnam in \citeg{AndPut} and then used by Bellissard and Savinien in \citeg{BelSav}.\\
Then we present the PV cohomology $H^*(\B_0; C(\Xi_\Delta,\ZZ))$ of pinwheel tilings which link the top \v{C}ech cohomology of $\Omega/S^1$ and the integer group of coinvariants of $\Xi_\Delta^2$.\\

\subsection{The pinwheel prototile space} \label{espacedesprototuiles}

\noindent
Let $\TT$ denote a pinwheel tiling (the one constructed in the beginning of this paper for example).\\
Let's remind some terminology :

\begin{defi} 
\begin{enumerate}
\item A \textbf{punctured} tile is an ordered pair consisting of a tile and a point in its interior.
\item A \textbf{prototile} of a tiling is an equivalence class of tiles (including the punctuation) up to direct isometries.
\item The \textbf{first corona} of a tile in a tiling $\TT$ is the union of the tiles of $\TT$ intersecting it.
\item A \textbf{collared prototile} of $\TT$ is the subclass of a prototile whose representatives have the same first corona up to direct isometries. 
\end{enumerate}
\end{defi}

\noindent
There is in fact $108$ collared prototiles for the pinwheel tiling (see \reffig{prototuilescouronnees} and \reffig{prototuilescouronnees2} on p.\pageref{prototuilescouronnees} and p.\pageref{prototuilescouronnees2}, where we have represented $54$ collared prototiles, the $54$ other collared prototiles are obtained by reflection).\\
We punctuate each prototile OF $\TT$ by the intersection point of the perpendicular bisector of the shortest side with the median from the vertex intersection of the hypotenuse and the shortest side:
\begin{center}
\begin{picture}(6,2)(0,0)
\linethickness{0.05cm}
\put(0.5,0.5){\line(1,0){2}}
\put(0.5,0.5){\line(0,1){1}}
\put(0.5,1.5){\line(2,-1){2}}
\put(0.5,1){\line(1,0){1.5}}
\put(0.5,1.5){\line(1,-1){1}}

\put(5.5,0.5){\line(-1,0){2}}
\put(5.5,0.5){\line(0,1){1}}
\put(3.5,0.5){\line(2,1){2}}
\put(5.5,1){\line(-1,0){1.5}}
\put(5.5,1.5){\line(-1,-1){1}}

\linethickness{0.1cm}
\put(1,1){\circle{0.05}}
\put(5,1){\circle{0.05}}
\end{picture}
\end{center}
If $\hat{t}$ is a prototile then $t$ will denote its representative that has its punctuation at the origin $O$.\\
We then take the simplicial structure defined in the next figure :
\begin{center}
\begin{picture}(6,4)(0,0)
\linethickness{0.05cm}
\put(1,1){\line(1,0){4}}
\put(1,1){\line(0,1){2}}
\put(1,3){\line(2,-1){4}}
\linethickness{0.01cm}
\put(1,3){\line(1,-2){1}}
\put(2,1){\line(0,1){1.5}}
\put(2,2.5){\line(2,-3){1}}
\put(2,1){\line(-1,1){1}}
\put(3,1){\line(0,1){1}}
\put(4,1){\line(-1,1){1}}
\put(4,1){\line(0,1){0.5}}
\end{picture}
\begin{picture}(6,4)(0,0)
\linethickness{0.05cm}
\put(5,1){\line(-1,0){4}}
\put(5,1){\line(0,1){2}}
\put(5,3){\line(-2,-1){4}}
\linethickness{0.01cm}
\put(5,3){\line(-1,-2){1}}
\put(4,1){\line(0,1){1.5}}
\put(3,1){\line(2,3){1}}
\put(4,1){\line(1,1){1}}
\put(3,1){\line(0,1){1}}
\put(2,1){\line(1,1){1}}
\put(2,1){\line(0,1){0.5}}
\end{picture}
\end{center}

\noindent
We can now build a finite CW-complex $\B_0(\TT)$, called prototile space, out of the collared prototiles by gluing them along their boundaries according to all the local configurations of their representatives in $\TT$ :

\begin{defi} \label{esp-proto}
Let $\hat{t}^c_j$, $j=1, \dots, N$, be the collared prototiles of $\mathcal{T}$ and let $t_j^c$ denote the representative of $\hat{t}_j^c$ that has its punctuation at the origin and the prototile orientation fixed during the construction of the canonical transversal.\\
The \textbf{collared prototile space} (or just the \textbf{prototile space}) of $\mathcal{T}$, $\mathcal{B}^c_0(\mathcal{T})$, is the quotient $CW$-complex
$$ \mathcal{B}^c_0(\mathcal{T}) = \coprod_{j=1}^N t_j^c / \sim,$$
where two $n$-cells $e^n_i \in t_i^c$ and $e^n_j \in t_j^c$ are identified if there exist direct isometries $(x_i,\theta_i)$,$(x_j,\theta_j) \in \RR^2 \rtimes SO(2)$ for which  $t^c_i.(x_i,\theta_i)$ and $t^c_j.(x_j,\theta_j)$ are tiles of $\mathcal{T}$ such that $e^n_i.(x_i,\theta_i)$ and $e^n_j.(x_j,\theta_j)$ coincide on the intersection of their $n$-skeletons.
\end{defi}

\noindent
The images of the tiles $t^c_j$ in $\B_0(\TT)$ will be denoted $\tau_j$ and still be called tiles.\\
We then have a projection map from $\Omega/S^1$ onto $\BBB$:

\begin{prop} \citeg{BelSav} \label{projec} There is a continuous map $\ppp:\Omega/S^1 \longrightarrow \BBB$ from the continuous hull quotiented by $S^1$ onto the collared prototile space.
\end{prop}

\noindent
\textbf{Proof :}
\begin{itemize}
\item[] Let ${\displaystyle \lambdo:\coprod_{j=1}^N t_j^c \rightarrow \BBB}$ be the quotient map and let ${\displaystyle \rho_0^c : \Omega/S^1 \rightarrow \coprod_{j=1}^N t_j^c}$ be defined as follows : take $[\omega] \in \Omega/S^1$ and $\omega_1$ a representative of this class.	If the origin $O$ belongs to the intersection of $k$ tiles $t^{\alpha_1}, \dots, t^{\alpha_k}$, in $\omega_1$, with $t^{\alpha_l} = t_{j_l}^c.\big(x_{\alpha_l}(\omega_1), \theta_{\alpha_l}(\omega_1) \big)$, $l=1, \dots, k$, then we set
$${\displaystyle \rho_0^c \big( [\omega] \big) = x_{\alpha_s}(\omega_1)}$$
which is in $t^c_{j_s}$ with $s=Min\{j_l : l=1, \ldots, k\}$.\\
This function depends on a choice of indice but this choice will vanish when we will look at the image in the quotient.\\
We remind that the $\EE^2$-action is given by : $\omega.(x,\theta) := R_{-\theta}(\omega - x)$ where $R_{-\theta}$ is the rotation in $\RR^2$ of angle $-\theta$ around the origin.\\
Thus, the definition of $\rho_0^c$ doesn't depend on the particular representative $\omega_1$ of $[\omega]$ chosen.\\
The map $\rho_0^c$ sends the origin of $\RR^2$, that lies in some tiles of a representative of $[\omega]$, to one of the corresponding tiles $t^c_j$'s at the corresponding position.\\
The projection $\ppp$ is then defined by :
$$\ppp : \begin{array}{ccc}
\Omega / S^1 & \longrightarrow & \BBB \\
 \text{} [\omega] & \longmapsto & \lambdo \circ \rho_0^c([\omega])\\
\end{array}.$$
We then have that, like in \citeg{BelSav}, $\ppp$ is well defined and continuous, noting that $[\omega']$ is in a neighborhood of $[\omega]$ in $\Omega/S^1$ if the representatives of $[\omega']$ and $[\omega]$ coincide on a big ball up to a small translation and up to rotations.
\end{itemize}

\begin{flushright}
$\square$
\end{flushright}

\noindent
For simplicity, the prototile space $\BBB$ is written $\BB$ and the projection $\ppp$ is written $\p_0$.\\
We denote $\Xi(\tau_j)$ the lift of the punctuation of $\tau_j$. This is a subset of the canonical transversal called the \textbf{acceptance zone} of the prototile $\hat{t}^c_j$. This subset contains all the tilings with the punctuation of a representative of $\hat{t}^c_j$ at the origin.\\
The $\Xi(\tau_j)$'s for $j = 1, \ldots, N_0$, form a clopen partition of the canonical transversal and are Cantor sets like $\Xi$.\\

\subsection{$\Omega/S^1$ as an inverse limit of supertile spaces}

We follow the guideline of \citeg{BelSav} to see our space $\Omega/S^1$ as the inverse limit of supertile spaces (see also \citeg{OrmRadSad} or \citeg{AndPut} to see the hull as an inverse limit).\\
Let $\TT$ be the pinwheel tiling constructed in the first section and set our simplicial decomposition on the prototiles.\\
We are going to define new finite $CW$-complexes $\B_k$, called supertile space of level $k$, associated to the $k$-supertiles of $\TT$.\\
The spaces $\B_k$ are built from the collared prototiles of an appropriate subtiling of $\TT$ , written $\TT_k$ below, in the same way as $\BB$ was built from the prototiles of $\TT$ in definition \refg{esp-proto}. The construction goes as follows.\\
As we said in the first section, $\TT$ can be decomposed uniquely in supertiles of level $k$, obtaining a repetitive non-periodic tiling $\TT_k$ of $\EE^2$-finite type and whose tiles are $k$-supertiles.\\
Each $k$-supertile is punctured by the punctuation of the tile in its middle as shown in the next figure:
\begin{center}
\begin{picture}(6,3.25)(0,0) \label{ponctuation}
\linethickness{0.05cm}
\put(1,1){\line(1,0){4}}
\put(1,1){\line(0,1){2}}
\put(1,3){\line(2,-1){4}}
\put(1,1){\line(1,2){0.8}}
\put(3,1){\line(-2,1){1.6}}
\put(3,1){\line(1,2){0.4}}
\put(3,1){\line(-3,4){1.2}}
\put(2,2){\circle{0.05}}
\end{picture}
\end{center}

\vspace*{-0.5cm}
\noindent
It's again the intersection point of the perpendicular bisector and the median of the appropriate side.\\
These $k$-supertiles are then compatible $CW$-complexes since they are made up of tiles of $\TT$ which are.

\begin{defi}
The \textbf{supertile space of level $\mathbf{k}$}, $\B_k$, is the collared prototile space of $\TT_k$ :
$$\B_k = \BB(\TT_k).$$
\end{defi}

\noindent
In fact, since all the tilings $\TT_k$ are the same, up to a dilatation, all the spaces $\B_k$ are homeomorphic but they will give us important informations on the cohomology of our space $\Omega/S^1$.\\
The images in $\B_k$ of the $k$-supertiles $p_j$ (tiles of $\TT_k$) are denoted $\pi_j$ and still called supertiles.\\
The projection $\p_{0,\TT_k}:\Omega/S^1 \longrightarrow \B_k$, built in proposition \refg{projec}, is denoted $\p_k$.\\
The map $F_k:\B_k \longrightarrow \B_0$ defined by $F_k:= \p_0 \circ \p_k^{-1}$ is well defined, onto and continuous (see \citeg{BelSav}). It projects $\B_k$ onto $\B_0$ in an obvious way : a point $x$ in $\B_k$ belongs to some supertile $\pi_j$, hence to some tile, and $F_k$ sends $x$ on the corresponding point in the corresponding tile $\tau_{j'}$.\\

\noindent
Let $p$ and $q$ be two integers such that $q \leqslant p$. The map $f_{q,p}:\B_p \rightarrow \B_q$ defined by $f_{q,p}:=F_q^{-1} \circ F_p=\p_q \circ \p_p^{-1}$ is well defined, onto and continuous.\\
As explained in \citeg{BelSav}, the family $(\B_p,f_{q,p})$ is a projective system.\\
To prove that $\Omega/S^1$ is the inverse limit $\displaystyle{ \lim_{\longleftarrow} (\B_p,f_{q,p})}$, we will use 
an important property of the pinwheel tiling : the $(l+1)$-supertiles have the same coronas as $l$-supertiles, up to a dilatation by a factor $\sqrt{5}$.\\
Hence, if $d_l$ denotes the distance between a $l$-supertile $p_l$ and the complementary of its first corona, the distance $d_{l+1}$ of the $(l+1)$-supertile which has the same first corona as the one of $p_l$ but dilated by $\sqrt{5}$, satisfies $d_{l+1}=\sqrt{5} d_l$ and thus, these distances goes to infinity. This point will allow us to prove the next theorem taking an appropriate sequence of supertile spaces.\\
We take the sequence of supertile spaces $\{\B_l,f_l\}_{l \in \NN}$ where, for $l \geqslant 1$, $\B_l$ is the space of $l$-supertiles and $f_l = f_{(l-1),l}$, with the convention that $f_0:=F_1$ and $\B_0$ is the prototile space.\\

\noindent
We then prove, as in \citeg{BelSav}, the following theorem:

\begin{theo} \label{limitprojective}
The inverse limit of the sequence $\{\B_l , f_l \}_{l \in \NN}$ is homeomorphic to the continuous hull of $\TT$ up to rotations :
$$\Omega / S^1 \cong \lim_{\longleftarrow} (\B_l,f_l) .$$
\end{theo}

\noindent
\textbf{Proof :}
\begin{itemize}
\item[] The homeomorphism is given by the map ${\displaystyle \p: \Omega / S^1 \longrightarrow  \lim_{\longleftarrow} (\B_l,f_l) }$, defined by  $\p([\omega])= \left ( \p_0([\omega]),\p_1([\omega]), \dots \right)$ with inverse  
$$\p^{-1}(x_0,x_1,\dots) = \cap\{\p_l^{-1}(x_l), l \in \NN\}.$$
The map $\p$ is surjective since each $\p_l$ is.\\
In fact, we have :
$$p_0^{-1}(x_0) \supset  p_1^{-1}(x_1) \supset p_2^{-1}(x_2) \supset \ldots$$
where each $p_i^{-1}(x_i)$ is a non empty compact subset of $\Omega/S^1$ and thus every tiling in the intersection above defines a lift of $(x_0,x_1,\dots)$.\\
\\
For the injectivity, let $\omega,\omega' \in \Omega$ be such that $\p([\omega])=\p([\omega'])$.\\
For each $l \in \NN$, $\p_l([\omega])=\p_l([\omega'])$ in some $l$-supertile $\pi_{l,j}$ of $\B_l$. This means that the two tilings agree, up to a rotation, on some translate of the supertile $p_{l,j}$ containing the origin.\\
Set then $r_l={\displaystyle \inf_{p \in \mathcal{P}_l^c} \inf_{x \in p} d_{\RR^2}(x,\partial C^1(p))}$ where $\mathcal{P}_l^c$ is the set of the collared $l$-supertiles of $\TT$ and $\partial  C^1(p)$ is the boundary of the first corona of $p$.\\
Since our tiling is of finite type, $r_l>0$ for all $l$.\\
Moreover, the two tilings agree, up to rotations, on the ball $B(0_{\RR^2},r_l)$ since $\B_l$ was built out of collared supertiles.\\
As mentioned earlier, $r_{l+1} = \sqrt{5}r_l$ thus $r_{l+1} = \sqrt{5}^{l+1}r_0$ and if we choose $l$ big enough, the two tilings agree on arbitrary large balls and so we have proved the injectivity of $\p$.\\
\\
$\p$ is then a bijection trivially continuous and since $ {\displaystyle \lim_{\longleftarrow} (\B_l,f_l)}$ is Hausdorff, $\p$ is in fact a homeomorphism.\\
\end{itemize}

\begin{flushright}
$\square$
\end{flushright}

\subsection{The PV cohomology}	

\noindent
We now turn to our first goal which was to compute the top \v{C}ech cohomology of $\Omega/S^1$ with integers coefficients and to prove that this was in fact the integer group of coinvariants $C(\Xi_\Delta^2,\ZZ)/\sim$ of a transversal under the partial action of the groupoid $\Omega \rtimes \RR^2 \rtimes S^1$ restricted to this transversal.\\
To prove this result, we will show that this cohomology is isomorphic to the PV cohomology introduced in \citeg{BelSav}, modified in order to take into account rotations.\\
We thus follow the approach of Bellissard and Savinien, introducing first the notion of oriented simplicial complexes, then defining the PV cohomology of the pinwheel tiling and finally proving that this is exactly the integer \v{C}ech cohomology of $\Omega/S^1$, thanks to the inverse limit found in the previous section.\\
The interesting point for us in this cohomology is that its cochains are in fact directly taken to be the continuous functions (with integer values) on the transversal $\Xi_\Delta^2$. We will then prove that the quotient of these cochains under the image of the differential of the PV cohomology is precisely the integer group of coinvariants on $\Xi_\Delta^2$.

\subsubsection{Oriented simplicial complexes and PV cohomology}

Here again, we will follow the presentation of Bellissard and Savinien in \citeg{BelSav} modifying the notions and the proofs to our oriented simplicial complexes (see \citeg{Hoc}).\\

\noindent
Given $n+1$ points $v_0, \dots, v_n$ in $\RR^m$ $m>n$, which are not colinear, let $[v_0,\dots,v_n]$ denote the $n$-simplex with vertices $v_0,\dots,v_n$.\\
Let $\Delta^n$ be the standard $n$-simplex :
$$\Delta^n = \left \{(x_0,\dots,x_n) \in \RR^{n+1} : \sum_{i=0}^n x_i = 1 \text{ and } x_i \geqslant 0 \text { for all } i \right \},$$
with vertices the unit vectors along the coordinate axis.\\
If one of the $n+1$ vertices of an $n$-simplex $[v_0,\dots,v_n]$ is deleted, the $n$ remaining vertices span a $(n-1)$-simplex, called a \textbf{face} of $[v_0,\dots,v_n]$.\\

\noindent
The \textit{boundary} of $\Delta^n$ is then the union of all the faces of $\Delta^n$ and is denoted $\partial \Delta^n$.\\
The \textit{interior} of $\Delta^n$ is then the open simplex $\overset{\!\!\!\! \circ}{\Delta^n}=\Delta^n \setminus \partial \Delta^n$.\\
A space $X$ is a \textbf{simplicial complex} if there is a collection of maps $\sigma_\alpha:\Delta^k \rightarrow X$, where $k$ depends on the index  $\alpha$, such that :
\begin{enumerate}[(i)]
	\item The restriction $\sigma_\alpha\text{}_{\mid \Delta^n}$ is injective. 

	\item Each restriction of $\sigma_\alpha$ to a face of $\Delta^n$ is one of the maps $\sigma_\beta:\Delta^{n-1} \rightarrow X$.

	\item For each $\alpha$ and $\beta$, $F_{\alpha,\beta}:=\sigma_\alpha(\Delta^n) \cap \sigma_\beta(\Delta^p)$ is a face of the two simplices $\sigma_\alpha(\Delta^n)$ and $\sigma_\beta(\Delta^p)$ and there is an affine map $l: \sigma_\alpha^{-1} (F_{\alpha,\beta}) \rightarrow \sigma_\beta^{-1}(F_{\alpha,\beta})$ such that $\sigma_\alpha\text{}_{\mid \sigma_\alpha^{-1} (F_{\alpha,\beta})} = \sigma_\beta\text{}_{\mid \sigma_\beta^{-1}(F_{\alpha,\beta})} \circ l$.
	\item A set $A \subset X$ is open iff $\sigma_\alpha^{-1}(A)$ is open in $\Delta^n$ for each $\sigma_\alpha$.
\end{enumerate}

\noindent
$\sigma_\alpha(\inte{\Delta^n})$ is called a \textbf{$\mathbf{n}$-cell} of the complex.\\

\noindent
We then obtain an \textbf{oriented simplex} from a $n$-simplex $\sigma=[v_0, \ldots, v_n]$ as follows : let fix an arbitrary ordering of the vertices $v_0, \ldots, v_n$. The equivalence class of even permutations of this fixed ordering is the \textbf{positively oriented simplex}, which we denote $+\sigma$. The equivalence class of odd permutations of the chosen ordering is the \textbf{negatively oriented simplex}, $- \sigma$.\\
An \textbf{oriented simplicial complex} is obtained from a simplicial complex by choosing an arbitrary fixed orientation for each simplex in the complex (this may be done without considering how the individual simplices are joined or whether one simplex is a face of another).\\
\\
To define an oriented simplicial structure on $\B_0$, we decompose each tile of the pinwheel tiling as follows :

\begin{center}
\begin{picture}(6,4)(0,0)
\linethickness{0.05cm}
\put(1,1){\line(1,0){4}}
\put(1,1){\vector(0,1){0.7}}
\put(1,1){\line(0,1){2}}
\put(1,1){\vector(1,0){0.7}}
\put(3,1){\vector(-1,0){0.8}}
\put(1,3){\vector(0,-1){0.8}}
\put(5,1){\vector(-1,0){0.9}}
\put(3,1){\vector(1,0){0.7}}
\put(1,3){\line(2,-1){4}}
\put(1,3){\vector(2,-1){0.7}}
\put(3,2){\vector(-2,1){0.6}}
\put(5,1){\vector(-2,1){0.7}}
\put(3,2){\vector(2,-1){0.8}}
\linethickness{0.01cm}
\put(1,3){\line(1,-2){1}}
\put(2,1){\line(0,1){1.5}}
\put(2,2.5){\line(2,-3){1}}
\put(2,1){\line(-1,1){1}}
\put(3,1){\line(0,1){1}}
\put(4,1){\line(-1,1){1}}
\put(4,1){\line(0,1){0.5}}
\end{picture}
\begin{picture}(6,4)(0,0)
\linethickness{0.05cm}
\put(5,1){\line(-1,0){4}}
\put(5,1){\vector(0,1){0.7}}
\put(5,1){\line(0,1){2}}
\put(5,1){\vector(-1,0){0.7}}
\put(3,1){\vector(1,0){0.8}}
\put(5,3){\vector(0,-1){0.8}}
\put(1,1){\vector(1,0){0.9}}
\put(3,1){\vector(-1,0){0.7}}
\put(5,3){\line(-2,-1){4}}
\put(5,3){\vector(-2,-1){0.7}}
\put(3,2){\vector(2,1){0.6}}
\put(1,1){\vector(2,1){0.7}}
\put(3,2){\vector(-2,-1){0.8}}
\linethickness{0.01cm}
\put(5,3){\line(-1,-2){1}}
\put(4,1){\line(0,1){1.5}}
\put(3,1){\line(2,3){1}}
\put(4,1){\line(1,1){1}}
\put(3,1){\line(0,1){1}}
\put(2,1){\line(1,1){1}}
\put(2,1){\line(0,1){0.5}}
\end{picture}
\end{center}
we take the orientation of $\RR^2$ for the $2$-cells and any orientation for the edges in the interior of our tiles.\\
We can thus see $\TT$ as an oriented simplicial decomposition of $\RR^2$.\\
We puncture each cell of each tile by the image under $\sigma_\alpha$ of the barycenter of $\Delta^n$ :

\begin{center}
\begin{picture}(6,4)(0,0)

\put(1,1){\line(1,0){4}}
\put(1,1){\line(0,1){2}}
\put(1,3){\line(2,-1){4}}
\put(1,3){\line(1,-2){1}}
\put(2,1){\line(0,1){1.5}}
\put(2,2.5){\line(2,-3){1}}
\put(2,1){\line(-1,1){1}}
\put(3,1){\line(0,1){1}}
\put(4,1){\line(-1,1){1}}
\put(4,1){\line(0,1){0.5}}

\linethickness{0.04cm}
\put(1.33,1.33){\circle{0.04}}
\put(1.33,2){\circle{0.04}}
\put(1.67,2.17){\circle{0.04}}
\put(2.33,1.5){\circle{0.04}}
\put(2.67,1.83){\circle{0.04}}
\put(3.33,1.33){\circle{0.04}}
\put(3.67,1.5){\circle{0.04}}
\put(4.33,1.17){\circle{0.04}}

\linethickness{0.06cm}
\put(1,1.5){\circle{0.05}}
\put(1,2.5){\circle{0.05}}
\put(1.5,1){\circle{0.05}}
\put(2.5,1){\circle{0.05}}
\put(3.5,1){\circle{0.05}}
\put(4.5,1){\circle{0.05}}
\put(1.5,1.5){\circle{0.05}}
\put(1.5,2){\circle{0.05}}
\put(2,1.75){\circle{0.05}}
\put(3,1.5){\circle{0.05}}
\put(3.5,1.5){\circle{0.05}}
\put(1.5,2.75){\circle{0.05}}
\put(2.5,2.25){\circle{0.05}}
\put(2.5,1.75){\circle{0.05}}

\linethickness{0.1cm}
\put(4,1.5){\circle{0.05}}
\put(2,2.5){\circle{0.05}}
\put(1,1){\circle{0.05}}
\put(1,2){\circle{0.05}}
\put(1,3){\circle{0.05}}
\put(2,1){\circle{0.05}}
\put(3,1){\circle{0.05}}
\put(4,1){\circle{0.05}}
\put(5,1){\circle{0.05}}
\put(3,2){\circle{0.05}}
\end{picture}
\begin{picture}(6,4)(0,0)
\put(5,1){\line(-1,0){4}}
\put(5,1){\line(0,1){2}}
\put(5,3){\line(-2,-1){4}}
\put(5,3){\line(-1,-2){1}}
\put(4,1){\line(0,1){1.5}}
\put(3,1){\line(2,3){1}}
\put(4,1){\line(1,1){1}}
\put(3,1){\line(0,1){1}}
\put(2,1){\line(1,1){1}}
\put(2,1){\line(0,1){0.5}}

\linethickness{0.04cm}
\put(4.67,1.33){\circle{0.04}}
\put(4.67,2){\circle{0.04}}
\put(2.67,1.33){\circle{0.04}}
\put(4.33,2.17){\circle{0.04}}
\put(3.67,1.5){\circle{0.04}}
\put(3.33,1.83){\circle{0.04}}
\put(2.33,1.5){\circle{0.04}}
\put(1.67,1.17){\circle{0.04}}

\linethickness{0.06cm}
\put(5,1.5){\circle{0.05}}
\put(5,2.5){\circle{0.05}}
\put(4.5,1){\circle{0.05}}
\put(3.5,1){\circle{0.05}}
\put(2.5,1){\circle{0.05}}
\put(1.5,1){\circle{0.05}}
\put(4.5,1.5){\circle{0.05}}
\put(4.5,2){\circle{0.05}}
\put(4,1.75){\circle{0.05}}
\put(3,1.5){\circle{0.05}}
\put(2.5,1.5){\circle{0.05}}
\put(4.5,2.75){\circle{0.05}}
\put(3.5,2.25){\circle{0.05}}
\put(3.5,1.75){\circle{0.05}}
\put(2.5,1.75){\circle{0.05}}
\put(1.5,1.25){\circle{0.05}}

\linethickness{0.1cm}
\put(2,1.5){\circle{0.05}}
\put(4,2.5){\circle{0.05}}
\put(5,1){\circle{0.05}}
\put(5,2){\circle{0.05}}
\put(5,3){\circle{0.05}}
\put(4,1){\circle{0.05}}
\put(3,1){\circle{0.05}}
\put(2,1){\circle{0.05}}
\put(1,1){\circle{0.05}}
\put(3,2){\circle{0.05}}
\end{picture}
\end{center}

\noindent
$\B_0$ is then a finite oriented simplicial complex and the maps $\sigma_\alpha : \Delta^n \rightarrow \B_0$ are the \textit{characteristic maps} of the $n$-simplices on $\B_0$.\\
We next define as in \citeg{BelSav} a new transversal (not immediatly related to the canonical transversal, see \refg{transversales}) which is crucial to define the PV cohomology:

\begin{defi}
The \textit{$\Delta$-transversal}, written $\Xi_\Delta$, is the subset of $\Omega/S^1$ formed by classes of tilings containing the origin on the punctuation of one of their cells.
\end{defi}

\noindent
The $\Delta$-transversal is the lift of the punctuation of the cells on $\B_0$.\\
It is partitioned by the lift of the punctuation of the $n$-cells, written $\Xi_{\Delta}^n$, i.e the subset of $\Omega/S^1$ consisting of classes of tilings containing the origin on the punctuation of one of their $n$-cells.\\
The $\Delta$-transversal is then a Cantor set (as the canonical transversal), and the $\Xi^n_\Delta$'s is a partition of it in clopen subsets.\\

\noindent
Let $\sigma$ be the characteristic map of a $n$-simplex $e$ on $\B_0$, denote $\Xi_\Delta(\sigma)$ the lift of the punctuation of $e$, and $\chi_\sigma$ its characteristic function in $\Xi_\Delta$ (i.e $\chi_\sigma([\omega])=1$ if and only if $\p_0([\omega]) = punct(e)$).\\
The subset $\Xi_\Delta(\sigma)$ is the \textit{acceptance zone} of $\sigma$.\\
Since $\Xi_\Delta(\sigma)$ is a clopen set, $\chi_\sigma \in C(\Xi_\Delta^n,\ZZ) \subset C(\Xi_\Delta,\ZZ)$.\\

\noindent
Consider $\sigma:\Delta^n \rightarrow \B_0$ the characteristic map of a $n$-simplex $e$ on $\B_0$, and let $\tau$ be a face of $\sigma$ with image cell $f$ in $\B_0$ (a face of $e$).\\
The simplices $e$ and $f$ on $\B_0$ are contained in some tile $\tau_j$. If we look at these simplices $e$ and $f$ as subsets of the tile $t_j$ in $\RR^2$, we can define a vector $x_{\sigma\tau}$ joining the punctuation of $f$ to the one of $e$.\\
The main issue with this construction is that, in the case of the pinwheel tiling, $x_{\sigma \tau}$ depends on the tile $\tau_j$ chosen if $e$ is a $1$-simplex (this vector is unique up to rotations). In the case of $2$-simplices, this vector is in fact unique since the quotient defining $\B_0$ only concerns the edges of the cells.\\
We must then find a way to choose such vectors in the case of $1$-simplices. There are many ways to do such a choice.\\
Here is one of them : if $e$ is a $1$-simplex, the vector is obtained by orienting the edge horizontally and from left to right in any tiling of $\Xi_\Delta(\sigma)$ :
\begin{center}
	\begin{picture}(1,1)(0,0)
			\put(0.5,0.5){\vector(1,0){0.6}}
			\put(0.5,0.5){\line(1,0){1}}
	\linethickness{0.06cm}
\put(0.5,0.5){\circle{0.05}}
\put(1.5,0.5){\circle{0.05}}

	\end{picture}
\end{center}

\noindent
Built in this manner, the vector $x_{\sigma\tau}$ are uniquely attached to $\sigma$ and $\tau$ as if we had a tiling where we only consider translations (as in \citeg{BelSav}).\\
We next define the "action" of $x_{\sigma \tau}$ on a function in $C(\Xi_\Delta,\ZZ)$ as follows : 
\begin{itemize}
\item if $\sigma$ is the characteristic map of a $2$-simplex $e$, $\tau$ a face of $\sigma$ and $f$ a function in $C(\Xi_\Delta(\tau),\ZZ)$ then, for each $[\omega] \in \Xi_\Delta(\sigma)$, we set
$$T^{x_{\sigma \tau}}f([\omega])=f([\omega_0+x_{\sigma \tau}]),$$
where $\omega_0 \in [\omega]$ is "well oriented", i.e the origin of the euclidean plane $E_2$ belongs to the punctuation of a cell $e$ (a triangle) of $\omega$ included in a unique tile (which is a direct isometry of the tile $t_j$ used to construct $x_{\sigma \tau}$); $\omega_0$ is then the tiling obtained by rotating $\omega$ to put this tile in this orientation (i.e such that this tile is a translation of $t_j$):
\begin{center}
\begin{picture}(6,2)(0,0)
\linethickness{0.05cm}
\put(0.5,0.5){\line(1,0){2}}
\put(0.5,0.5){\line(0,1){1}}
\put(0.5,1.5){\line(2,-1){2}}
\put(5.5,0.5){\line(-1,0){2}}
\put(5.5,0.5){\line(0,1){1}}
\put(3.5,0.5){\line(2,1){2}}
\end{picture}
\end{center}

			\item If $\sigma$ is the characteristic map of a $1$-simplex, $\tau$ a face of $\sigma$ and $f$ a function in $C(\Xi_\Delta(\tau),\ZZ)$ then, for each $[\omega] \in \Xi_\Delta(\sigma)$, we put
$$T^{x_{\sigma \tau}}f([\omega])=f([\omega_0+x_{\sigma \tau}]),$$
where $\omega_0 \in [\omega]$ is "well oriented", i.e the origin of $E_2$ belongs to the punctuation of an edge of $\omega$; $\omega_0$ is then obtained by rotating $\omega$ to put this edge horizontally and the orientation in the positive direction (from left to right).\\
\end{itemize} 

\noindent
We then define important operators, useful to define the differential of the PV cohomology :

\begin{defi} \label{opetheta}
Let $\sigma$ and $\tau$ be the characteristic map of a $n$-simplex, resp. a $k$-simplex, on $\B_0$.\\
We define the operator $\theta_{\sigma \tau}$ on $C(\Xi_\Delta,\ZZ)$ by :
$$\theta_{\sigma \tau} = \left\{ \begin{array}{cl}
\chi_\sigma T^{x_{\sigma \tau}}\chi_\tau & \text{ if } \tau \subset \partial \sigma \text{ and } n=1, 2 \\
0 & \text{ else }
\end{array} \right.$$
where $\tau \subset \partial \sigma$ means that $\tau$ is a face of $\sigma$ of codimension $1$.
\end{defi}

\noindent
This operator is easy to describe : if $f \in C(\Xi_\Delta(\tau),\ZZ)$ and $[\omega] \in \Xi_\Delta$ then $\theta_{\sigma \tau}(f)([\omega])$ is $0$ if $\omega$ is not in $\Xi_\Delta(\sigma)$ or if $\tau$ is not a face of $\sigma$ of codimension $1$ or if $n \neq 1, 2$, and else,  $\theta_{\sigma \tau}(f)([\omega])$  is equal to the value of $f$ on the class of the tiling obtained from $\omega$ by a direct isometry sending the origin to the punctuation of the face $\tau$ of $\sigma$ and rotating $\tau$ such that it has the "good" orientation.\\

\noindent
We can then define the PV cohomology of pinwheel tilings.\\
Let $\mathcal{S}^n_0$ denote the set of the characteristic maps $\sigma:\Delta^n \rightarrow \B_0$ of $n$-simplices on $\B_0$ and $\mathcal{S}_0$ the union of the $\mathcal{S}_0^n$'s. The group of simplicial $n$-chains on $\B_0$, $C_{0,n}$, is the free abelian group with basis $\mathcal{S}_0^n$.\\

\noindent
Before defining the PV cohomology, we need one more definition, the incidence number:

\begin{defi}
Let $\sigma$ and $\tau$ be two simplices of dimension $n$ and $n-1$ respectively, the \textbf{incidence number} $[\sigma,\tau]$ is defined by :
$$\begin{array}{ccccc}
 & [\sigma,\tau]  & = & \pm 1 &  \text{ if } \tau \subset \partial \sigma \\
 & [\sigma,\tau]  & = & 0 &  \text{ else } \\
\end{array}$$
If $\tau \subset \partial \sigma$, then  $[\sigma,\tau]$ is $1$ if $\tau$ is a \textit{positively oriented face} of $\sigma$ i.e the orientation of $\tau$ coincides with that induced by $\sigma$ on its face $\tau$ and this number is $-1$ if $\tau$ is a \textit{negatively oriented face} of $\sigma$.
\end{defi}

\begin{defi}
The \textit{PV cohomology} of $\Omega / S^1$ is the cohomology of the differential complex $\{C_{PV}^n,d_{PV}^n\}$, with :
\begin{enumerate}
		\item the PV cochain groups are the groups of continuous integer valued functions on $\Xi_\Delta^n$ : $C^n_{PV}=C(\Xi_\Delta^n,\ZZ)$ for $n=0,1,2$,
		\item the PV differential, $d_{PV}$, is defined by the sum over $n=1,2$, of the operators :
$$d^n_{PV}:\left\{ \begin{array}{c}
C_{PV}^{n-1}  \longrightarrow  C^n_{PV}\\
d^n_{PV}={\displaystyle \sum_{\sigma \in \mathcal{S}^n_0} \sum_{i=0}^n [\sigma,\partial_i \sigma] \theta_{\sigma\partial_i\sigma}}
\end{array}\right. .$$
\end{enumerate}
\end{defi}

\noindent
The "simplicial form" of $d^n_{PV}$ easily implies $d^2_{PV} \circ d^1_{PV}=0$.\\
This comes from the fact that, for each simplex $\sigma$ (see \citeg{Hoc}) :
$$\sum_{i,j} [\sigma,\partial_i\sigma][\partial_i \sigma,\partial_j \partial_i \sigma] = 0.$$
We shall also call this cohomology the \textbf{PV cohomology of $\TT$}.\\
We denote it $H^*_{PV} \big ( \B_0;C(\Xi_\Delta,\ZZ)\big)$.\\

\noindent
The next subsection will then prove the following theorem :

\begin{theo} \label{isocechpimsner}
The integer \v{C}ech cohomology of $\Omega/S^1$ is isomorphic to the PV cohomology of $\TT$ :
$$\check{H}^*(\Omega/S^1 ; \ZZ) \cong H^*_{PV} \big ( \B_0;C(\Xi_\Delta,\ZZ)\big).$$
\end{theo}

\subsubsection{Proof of theorem \refg{isocechpimsner}}

\noindent
Once again, we follow the guideline of the paper \citeg{BelSav}.\\
We define first a PV cohomology for the $\B_p$'s, written $H^*_{PV}(\B_0;C(\Sigma_p,\ZZ))$. We show that this cohomology is in fact the simplicial cohomology of $\B_p$ in the proposition \refg{pimsnersimpliciale} and then we prove that the PV cohomology of $\TT$ is isomorphic to the direct limit of the PV cohomologies of the supertile spaces sequence used in theorem \refg{limitprojective}.\\

\noindent
Denote $\mathcal{S}_p^n$ the set of all the characteristic maps $\sigma_p:\Delta^n \rightarrow \B_p$ of the $n$-simplices on $\B_p$, and $\mathcal{S}_p$ the union of the $\mathcal{S}^n_p$'s.\\
The group of simplicial $n$-chains on $\B_p$, $C_{p,n}$, is the free abelian group with basis $\mathcal{S}_p^n$.\\

\noindent
As above, if $\sigma_p$ is a simplex on $\B_p$, write $\Xi_{p,\Delta}(\sigma_p)$ for the lift of the punctuation of its image in $\B_p$ and $\chi_{\sigma_p}$ its characteristic map. $\Xi_{p,\Delta}(\sigma_p)$ is the \textit{acceptance zone} of $\sigma_p$. It's a clopen subset of the $\Delta$-transversal.

\begin{lem} \label{partition}
Given a simplex $\sigma$ on $\B_0$, its acceptance zone is partitioned by the acceptance zone of its preimages in $\B_p$ :
$$\Xi_\Delta(\sigma) = \bigsqcup_{\sigma_p \in F_{p\#}^{-1}(\sigma)} \Xi_{p,\Delta}(\sigma_p),$$
with $F_{p\#}:\mathcal{S}_p^n \rightarrow \mathcal{S}_0^n$ the map induced by $F_p$.
\end{lem}

\noindent
The proof is exactly the same that the one in \citeg{BelSav}.\\ 

\noindent
We note that the union over $\sigma \in \mathcal{S}^n_0$ of the $F^{-1}_{p\#}(\sigma)$'s is $\mathcal{S}^n_p$.\\

\noindent
We then denote $C^n_p$ the simplicial $n$-cochain group $Hom(C_{p,n},\ZZ)$, which is the dual of the simplicial $n$-chain group $C_{p,n}$.\\
We can represent it faithfully on the group of continuous function with integer values on the $\Delta$-transversal $C(\Xi_\Delta,\ZZ)$ by :
$$\rho_{p,n}: \left\{ \begin{array}{ccc}
C_p^n & \longrightarrow & C(\Xi_\Delta^n,\ZZ) \\
\psi & \longmapsto & {\displaystyle \sum_{\sigma_p \in \mathcal{S}^n_p}} \psi(\sigma_p)\chi_{\sigma_p}
\end{array} \right..$$
We denote $C(\Sigma^n_p,\ZZ)$ the image of this representation.\\
$\rho_{p,n}$ is an isomorphism on its image $C(\Sigma_p^n,\ZZ)$, its inverse is defined as follows: given $\phi={\displaystyle \sum_{\sigma_p \in \mathcal{S}_p^n}} \phi_{\sigma_p}\chi_{\sigma_p}$, where $\phi_{\sigma_p}$ is an integer, $\rho_{p,n}^{-1}(\phi)$ is the group homomorphism from $C_{p,n}$ to $\ZZ$ whose value on the basis simplex $\sigma_p$ is $\phi_{\sigma_p}$.\\

\noindent
Consider the characteristic map $\sigma_p$ of a $n$-simplex $e_p$ on $\B_p$. This simplex is contained in some supertile $\pi_j$. Viewing $e_p$ as a subset of the supertile $p_j$ in $\RR^2$, we can again define, similarly to the method used for the PV cohomology, the vector $x_{\sigma_p \partial_i\sigma_p}$, for $i=1,\dots,n$, joining the punctuation of the $i$-th face $\partial_i e_p$ to the punctuation of $e_p$. \\
Since $F_p$ preserves the orientation of the simplices, these vectors $x_{\sigma_p \partial_i \sigma_p}$ are identical, for each $\sigma_p$ in the preimage of the characteristic map $\sigma$ of a simplex $e$ on $\B_0$, and they are, in fact, equal to the vector $x_{\sigma \partial_i \sigma}$ used in the definition \refg{opetheta} of the operator $\theta_{\sigma \partial_i\sigma}$. In the same way, we define the operators $\theta_{\sigma_p \partial_i\sigma_p}$ as the operators $\chi_{\sigma_p} T^{x_{\sigma_p \partial_i \sigma_p}}\chi_{\partial_i\sigma_p}$.\\
Using the relation $T^{x_{\sigma \partial_i\sigma}} \chi_{\partial_i \sigma} = \chi_{\sigma}T^{x_{\sigma \partial_i \sigma}}$ and lemma \refg{partition}, we obtain :
$$\theta_{\sigma\partial_i\sigma} = \sum_{\sigma_p \in F_{p\#}^{-1}(\sigma)} \theta_{\sigma_p \partial_i \sigma_p}.$$
Hence, the PV differential can be written :
$$d^n_{PV}=\sum_{\sigma_p \in \mathcal{S}^n_p} \sum_{i=0}^{n} [\sigma_p,\partial_i \sigma_p] \theta_{\sigma_p \partial_i \sigma_p},$$
and this defines a differential from $C(\Sigma_p^{n-1},\ZZ)$ to $C(\Sigma_p^n,\ZZ)$ (since $F_p$ preserves the orientation of the simplices, for each $\sigma_p \in F_{p\#}^{-1}(\sigma)$, $[\sigma_p,\partial_i \sigma_p] = [\sigma, \partial_i\sigma]$, which justifies the definition of this differential).

\begin{defi}
Set $C^n_{PV}(p)=C(\Sigma_p^n,\ZZ)$, for $n=0,1,2$. The PV cohomology of the supertile space $\B_p$, written $H^*_{PV}\big( \B_0;C(\Sigma_p,\ZZ) \big)$, is the cohomology of the differential complex $\{C^n_{PV}(p),d^n_{PV}\}$.
\end{defi}

\noindent
As in \citeg{BelSav}, we obtain one of the two crucial propositions for the proof of \refg{isocechpimsner}:

\begin{prop} \label{pimsnersimpliciale}
The PV cohomology of the supertile space $\B_p$ is isomorphic to its integer simplicial cohomology :
$$H^*_{PV}\big( \B_0;C(\Sigma_p,\ZZ) \big) \cong H^*\big( \B_p;\ZZ \big).$$
\end{prop}

\noindent
\textbf{Proof :}
		\begin{enumerate}
			\item[] $\rho_{p,n}$ is an isomorphism from the simplicial cochain group $C_p^n$ onto the PV cochain group $C^n_{PV}(p)$.\\
Let $\phi$ be an element in $C^{n-1}_{PV}(p)$ and $\sigma_p$ a $n$-simplex in $\mathcal{S}_p^n$, the differential of $\phi$ is then given by :
$$d^n_{PV}\phi(\sigma_p) =  \sum_{i=0}^n [\sigma_p,\partial_i \sigma_p] \phi(\partial_i \sigma_p).$$
On the other hand, the simplicial differential of some $\psi \in C_p^{n-1}$ is :
$$\delta^n\psi(\sigma_p) = \sum_{i=1}^n [\sigma_p,\partial_i] \psi(\partial_i \sigma_p).$$
thus, we have $d^n_{PV} \circ \rho_{p,n-1} = \rho_{p,n} \circ \delta^n$ for $n=1,2$.\\
So the $\rho_{p,n}$'s give a chain map and thus induce isomorphisms $\rho_{p,n}^*$'s between the $n$-th cohomology groups.
		\end{enumerate}
\begin{flushright}
$\square$
\end{flushright}

\noindent 
The following lemma is trivial using theorem \refg{limitprojective} and is useful to prove the second crucial proposition for the proof of theorem \refg{isocechpimsner} :

\begin{lem} \label{limitdeltatransv}
Let $\{\B_l,f_l\}$ be the sequence used in theorem \refg{limitprojective}. We have :
$$\Xi_\Delta \cong \lim_{\longleftarrow}(\mathcal{S}_l,f_l)$$
and
$$C(\Xi_\Delta,\ZZ) \cong \lim_{\longrightarrow}\Big( C(\Sigma_l,\ZZ),f^l \Big),$$
where the $f^l$'s are the duals of the $f_l$'s.
\end{lem}

\begin{prop} \label{pimsnerlimitinductive}
Let $\{\B_l,f_l\}$ be as in the previous lemma. There is an isomorphism :
$$H^*_{PV}\big( \B_0;C(\Xi_\Delta,\ZZ) \big) \cong \lim_{\longrightarrow} \Big( H^*_{PV}\big( \B_0;C(\Sigma_p,\ZZ) \big ) , f_l^* \Big).$$
\end{prop}

\noindent
\textbf{Proof :}
		\begin{enumerate}
			\item[] By the previous lemma \refg{limitdeltatransv}, the cochain group $C^n_{PV}$ are the direct limits of the cochain groups $C^n_{PV}(l)$ of the supertile spaces $\B_l$.\\
Consider $f_l^\#:C^n_{PV}(l) \longrightarrow C^n_{PV}(l+1)$ the map induced by $f_l$ on the PV cochain groups.\\
Since the differential $d_{PV}$ is the same for the complexes of each supertile space $\B_l$, it's enough to check that the following diagram is commutative
$$\xymatrix{
\cdots \ar[r] & C_{PV}^{n-1}(l) \ar[r]^{d^n_{PV}}  \ar[d]^{f^\#_l}  & C^n_{PV}(l)  \ar[r]  \ar[d]^{f_l^\#} & \cdots \\
\cdots \ar[r] & C_{PV}^{n-1}(l+1) \ar[r]^{d^n_{PV}}  & C^n_{PV}(l+1) \ar[r] & \cdots 
}$$
and this is easy using the relations 
$$(f_l^\#\phi)(\partial_i \sigma_l) = \phi(f_{l \#}(\partial_i \sigma_l) ) \text{ and } f_{l \#}(\partial_i \sigma_l) = \partial_i ( f_{l \#}(\sigma_l)).$$
\end{enumerate}
\begin{flushright}
$\square$
\end{flushright}

\noindent
To end the proof of theorem \refg{isocechpimsner}, we then use the fact that \v{C}ech cohomology sends inverse limits on direct limits, that \v{C}ech cohomology is isomorphic to simplicial cohomology for $CW$-complexes (the $\B_l$'s) and that the simplicial cohomology of $\B_l$ is isomorphic to the PV cohomology of $\B_l$ by proposition \refg{pimsnersimpliciale}.\\
We then conclude that the direct limit of these groups is the PV cohomology of $\TT$ by proposition \refg{pimsnerlimitinductive}

\subsection{The integer group of coinvariants of the pinwheel tiling}

\noindent
In this section, we show that the top PV cohomology $H^2_{PV}(\B_0;C(\Xi_\Delta,\ZZ))$ is isomorphic to the integer group of coinvariants $C(\Xi_\Delta^2,\ZZ)/\sim$ of the transversal $\Xi_\Delta^2$ (definition follows) and thus, the top integer \v{C}ech cohomology of the hull of a pinwheel tiling is also isomorphic to this group which is a first step in the proof of the gap-labelling.\\

\noindent
In the same way we identified the canonical transversal $\Omega/S^1$ to a subset of $\Omega$, we can identify $\Xi^2_\Delta$ and $\Xi_\Delta^1$ with a subset of $\Omega$.\\
Indeed, it suffices to choose in each class $[\omega] \in \Xi^n_\Delta$ ($n=1,2$) a representative with the "good orientation" (the one used to define the vector $x_{\sigma \tau}$).\\
We also obtain a representation of the clopens $\Xi_{p,\Delta}$ in $\Omega$ in this way.\\

\noindent
We will denote in the same way these subsets of $\Omega/S^1$ and the corresponding subsets of $\Omega$. When we take a class in $\Xi_\Delta^n$ or  $\Xi_{p,\Delta}$, we implicitly mean that we are in $\Omega/S^1$ and if we take a tiling in $\Xi_\Delta^n$ or $\Xi_{p,\Delta}$, we consider the subset of $\Omega$.\\

\noindent
Let's begin by the definition of the integer group of coinvariants of the canonical transversal $\Xi$ and of the $n$-th $\Delta$-transversal $\Xi_\Delta^n$ ($n=1,2$).\\
We let $\Xi_0$ denote either $\Xi$ or $\Xi_\Delta^n$ ($n=1,2$) (subsets of $\Omega$).\\
Let $\omega_0 \in \Xi_0$ and $\mathcal{A}_{\omega_0}$ a patch of $\omega_0$ around the origin, we define 
$$U(\omega_0,\mathcal{A}_{\omega_0}):=\{\omega' \in \Xi_0 \mid \omega_0 \text{ and } \omega' \text{ coincide on } \mathcal{A}_{\omega_0}\},$$
which is a clopen subset of $\Xi_0$.\\
Viewing this subset in $\Omega/S^1$, this clopen becomes
$$V([\omega_0],\mathcal{A}_{\omega_0}) := \{[\omega_1] \in \Xi_0 \mid \omega_0 \text{ and } \omega_1 \text{ coincide on } \mathcal{A}_{\omega_0} \text{ up to rotations }\}.$$
By hypothesis, our tiling $\TT$ is of finite $\RR^2 \rtimes SO(2)$-type hence the family 
$$\{U(\omega_0,\mathcal{A}_{\omega_0}) , \omega_0 \in \Xi_0, \mathcal{A}_{\omega_0} \text{ patch of size } k \text{ of } \omega_0 \text{ around } O \, , k \in \NN \}$$
 is countable.\\
We define for each $\omega_0 \in \Xi_0$:
$$\mathcal{G} \Big( U(\omega_0,\mathcal{A}_{\omega_0}) \Big)=\{(x,\theta) \in \RR^2 \rtimes SO(2) \mid \omega_0.(x,\theta) \in \Xi_0 \text{ and } x \in \mathcal{A}_{\omega_0}\},$$
where $x \in \mathcal{A}_{\omega_0}$ means that $x$ is a vector in $\RR^2$ contained in the subset of $\RR^2$ defined by $\mathcal{A}_{\omega_0}$.\\
$\mathcal{G} \Big( U(\omega_0,\mathcal{A}_{\omega_0}) \Big)$ is defined like this because in fact, if we take a tiling $\omega'$ in the clopen set $U(\omega_0,\mathcal{A}_{\omega_0})$, we know this tiling only on the patch $\mathcal{A}_{\omega_0}$ and thus we have :
$$\forall (x,\theta) \in \mathcal{G} \Big( U(\omega_0,\mathcal{A}_{\omega_0}) \Big), \, \forall \, \omega' \in U(\omega_0,\mathcal{A}_{\omega_0}), \, \omega'.(x,\theta) \in \Xi_0.$$

\noindent
The integer group of coinvariants of $\Xi_0$ is then the quotient of $C(\Xi_0,\ZZ)$ by the subgroup $H_{\Xi_0}$ spanned by the family
$$\{\chi_{U(\omega_0,\mathcal{A}_{\omega_0})}-\chi_{U(\omega_0,\mathcal{A}_{\omega_0}).(x,\theta)} \mid (x,\theta) \in \mathcal{G}(U(\omega_0,\mathcal{A}_{\omega_0})), \omega_0 \in \Xi_0\}.$$

\begin{theo}
$$H^2_{PV}(\B_0;C(\Xi_\Delta,\ZZ)) \cong C(\Xi^2_\Delta,\ZZ)/ H_{\Xi_\Delta^2}$$
\end{theo}

\noindent
\textbf{Proof :}
		\begin{enumerate}
			\item[] By definition, this cohomology is $C(\Xi^2_\Delta,\ZZ)/ \text {Im}(d^2_{PV})$.\\
We prove that $\text{Im}(d^2_{PV}) = H_{\Xi_\Delta^2}$.\\

\item Consider $f \in \text{Im} (d^2_{PV})$, there is some $g \in C(\Xi_\Delta^1,\ZZ)$ such that $f=d^2_{PV}(g)$.\\
Since $C(\Xi_\Delta^1,\ZZ)$ is generated by characteristic functions of the form $\chi_{U(\omega,\mathcal{A}^\omega_\tau)}$ with $\tau \in \mathcal{S}^1_p$ for $p \in \NN$, $\omega \in \Xi_{p,\Delta}(\tau)$ and $\mathcal{A}^\omega_\tau$ a patch of $\omega$ around the origin (in fact around the edge which projects on $\tau$ in $\B_p$) large enough to cover the first corona(s) of the (two) supertile(s) surrounding $\tau$ in $\omega$, it is enough to prove the result for $f=d^2_{PV}(\chi_{U(\omega,\mathcal{A}^\omega_\tau)})$ with $\tau$ in $\mathcal{S}^1 := \bigcup \mathcal{S}^1_p$.\\
Let thus $\tau$ be some element in $\mathcal{S}^1_p$, $\omega$ a tiling in $\Xi_{p,\Delta}(\tau)$ and $\mathcal{A}^\omega_\tau$ a patch surrounding the origin in $\omega$ (and thus surrounding the "edge" $\tau$) large enough, then $\chi_{U(\omega,\mathcal{A}^\omega_\tau)}$ is the characteristic function of the set of all the tilings in $\Omega$ with the origin on the punctuation of $\tau$, $\tau$ having the "good" orientation and which coincide on the patch $\mathcal{A}^\omega_\tau$ with $\omega$.\\
We remark that if $\tau'$ is the characteristic map of another $1$-simplex in $\mathcal{S}^1_p$, then 
$$\chi_{\tau'} \chi_{U(\omega,\mathcal{A}^\omega_\tau)} = \delta_{\tau \tau'} \chi_{U(\omega,\mathcal{A}^\omega_\tau)}.$$
As $\mathcal{A}^\omega_\tau$ was chosen large enough, this patch characterize the (two) collared supertile(s) surrounding the simplex corresponding to $\tau$ in the tilings of $U(\omega,\mathcal{A}^\omega_\tau)$. Thus, denoting $\sigma_0$ and $\sigma_1$ the two characteristic maps of the $2$-simplices having $\tau$ as an edge in $\omega$ and the above collared supertile(s) (respectively) as supertiles in $\B_p$, we have : \\
for each $\sigma \in \mathcal{S}_p^2$
$$\chi_{\sigma} \chi_{U(\omega^0,\mathcal{A}^0_\tau)} = \delta_{\sigma \sigma_0}  \chi_{U(\omega^0,\mathcal{A}^0_\tau)}$$
and
$$\chi_{\sigma} \chi_{U(\omega^1,\mathcal{A}^1_\tau)} = \delta_{\sigma \sigma_1}  \chi_{U(\omega^1,\mathcal{A}^1_\tau)}$$
where $\omega^i=R_{\theta_i}(\omega)-x_{\sigma_i \tau}$ ($i=0,1$), $\mathcal{A}^i_\tau = R_{\theta_i}(\mathcal{A}^\omega_\tau)-x_{\sigma_i \tau}$ and $\omega^i$ is in $\Xi_{p,\Delta}(\sigma_i)$ ($\omega^i$ are the tilings in $\Xi_{p,\Delta}(\sigma_i)$ obtained from $\omega$ by a direct isometry taking the origin on the punctuation of $\sigma_i$ and rotating $\omega^i$ to obtain the good orientation).\\
With these three remarks, we easily see that :
$$d^2_{PV}(\chi_{U(\omega,\mathcal{A}^\omega_\tau)}) = \pm ( \chi_{U(\omega^0,\mathcal{A}^0_\tau)} - \chi_{U(\omega^1,\mathcal{A}^1_\tau)}),$$
with $\omega^1=\omega^0.(R_{\theta_0-\theta_1}(x_{\sigma_1 \tau})-x_{\sigma_0 \tau},\theta_0-\theta_1)$ and $\mathcal{A}^1_\tau$ of the same form.\\
Hence
$$d^2_{PV}(\chi_{U(\omega,\mathcal{A}^\omega_\tau)}) = \pm( \chi_{U(\omega^0,\mathcal{A}^0_\tau)} - \chi_{U(\omega^0,\mathcal{A}^0_\tau).(y_{\sigma_0 \sigma_1},\theta_{\sigma_0 \sigma_1)}}),$$
where $y_{\sigma_0 \sigma_1} = R_{\theta_0-\theta_1}(x_{\sigma_1 \tau})-x_{\sigma_0 \tau}$, $\theta_{\sigma_0 \sigma_1}= \theta_0-\theta_1$ and thus  $(y_{\sigma_0 \sigma_1},\theta_{\sigma_0 \sigma_1})$ is in $\mathcal{G} \Big( U(\omega_0,\mathcal{A}_{\omega_0}) \Big)$.\\
We thus have proved that $d^2_{PV}(\chi_{U(\omega,\mathcal{A}_\tau)}) \in H_{\Xi_\Delta^2}$ and thus the inclusion $\text{Im}(d^2_{PV}) \subset H_{\Xi_\Delta^2}$.\\

			\item The inclusion in the other direction is easy to obtain by reasoning on generators.\\
Let $\chi_{U(\omega,\mathcal{A}_\omega)}-\chi_{U(\omega,\mathcal{A}_\omega).(x,\theta)}$ be a generator of $H_{\Xi_\Delta^2}$.\\
Covering $\mathcal{A}_\omega$ by supertiles large enough, we can suppose that $\mathcal{A}_\omega$ is in fact a collared supertile or a union of two collared supertiles of same level or a star of collared supertiles of same level (a star of supertiles is all the supertiles surrounding a fixed vertex).\\

\noindent
Consider $\omega \in \Xi_\Delta^2$ and $\mathcal{A}_\omega$ a patch in $\omega$ consisting of $p$-supertiles around the origin of the above form, called \textit{collared patch}.\\
Thanks to the form of collared patches, we can find a sequence of tilings $\omega_0, \ldots, \omega_n$ such that $\omega_0=\omega$, $\omega_n=\omega.(x,\theta)$, $\omega_i = \omega.(x_i,\theta_i)$ ($i=1, \ldots, n$) with $(x_i,\theta_i) \in \mathcal{G} \Big( U(\omega,\mathcal{A}_\omega) \Big)$ and such that the $2$-simplex in $\omega_i$ containing the origin have a common edge with the $2$-simplex in $\omega_{i-1}$ which contained the origin (see the figure below).\\
If we take, for example, the patch below, with the points representing $\omega$ and $\omega.(x,\theta)$ (in fact, on this figure, we have represented $\omega$ and $\omega-x$, thus you must think that we are representing tilings by points and then you must rotate the tiling to put the tile containing the origin in the "good orientation" to obtain $\omega$ and $\omega.(x,\theta)$) :
\begin{center}
\hspace*{-2cm}\begin{picture}(10,5)(-1,-1) 
\linethickness{0.02cm}
\put(0,0){\line(1,0){10}}
\put(0,0){\line(2,1){8}}
\put(10,0){\line(-1,2){2}}
\put(4,0){\line(0,1){2}}
\put(8,0){\line(0,1){4}}
\put(4,2){\line(1,0){5}}
\put(8,0){\line(-2,1){4}}
\put(2,0){\line(0,1){1}}
\put(6,0){\line(0,1){3}}
\put(1,0){\line(1,1){1}}
\put(3,0){\line(1,2){1}}
\put(3,0){\line(1,1){1}}
\put(5,0){\line(-1,1){1}}
\put(5,0){\line(-1,2){1}}
\put(7,0){\line(-1,1){2}}
\put(9,0){\line(-1,1){2}}
\put(8,0){\line(-1,2){1}}
\put(10,0){\line(-2,1){2}}
\put(6,3){\line(-1,-1){1}}
\put(7,2){\line(1,2){1}}
\put(7,2){\line(1,1){1}}
\put(9,2){\line(-1,1){1}}
\put(1,0){\line(0,1){0.5}}
\put(3,0){\line(0,1){1.5}}
\put(5,0){\line(0,1){2.5}}
\put(7,0){\line(0,1){3.5}}
\put(8,1){\line(1,0){1.5}}
\put(2,0){\line(2,3){1}}
\put(6,0){\line(-2,3){1}}
\put(6,2){\line(2,-3){1}}
\put(6,2){\line(2,3){1}}
\put(8,2){\line(3,-2){1.5}}
\put(8,3){\line(1,0){0.5}}

\linethickness{0.1cm}
\put(2.67,0.5){\circle{0.04}}
\put(6.67,2.5){\circle{0.04}}

\linethickness{0.07cm}
\put(2.67,0.5){\vector(2,1){4}}
\end{picture}
\end{center}

\noindent
We may thus represent the sequence of tilings "joining" $\omega$ and $\omega.(x,\theta)$ by the following sequence of points (there isn't a unique path) :

\begin{center}
\hspace*{-2cm} \begin{picture}(10,5)(-1,-1) 
\linethickness{0.02cm}
\linethickness{0.02cm}
\put(0,0){\line(1,0){10}}
\put(0,0){\line(2,1){8}}
\put(10,0){\line(-1,2){2}}
\put(4,0){\line(0,1){2}}
\put(8,0){\line(0,1){4}}
\put(4,2){\line(1,0){5}}
\put(8,0){\line(-2,1){4}}
\put(2,0){\line(0,1){1}}
\put(6,0){\line(0,1){3}}
\put(1,0){\line(1,1){1}}
\put(3,0){\line(1,2){1}}
\put(3,0){\line(1,1){1}}
\put(5,0){\line(-1,1){1}}
\put(5,0){\line(-1,2){1}}
\put(7,0){\line(-1,1){2}}
\put(9,0){\line(-1,1){2}}
\put(8,0){\line(-1,2){1}}
\put(10,0){\line(-2,1){2}}
\put(6,3){\line(-1,-1){1}}
\put(7,2){\line(1,2){1}}
\put(7,2){\line(1,1){1}}
\put(9,2){\line(-1,1){1}}
\put(1,0){\line(0,1){0.5}}
\put(3,0){\line(0,1){1.5}}
\put(5,0){\line(0,1){2.5}}
\put(7,0){\line(0,1){3.5}}
\put(8,1){\line(1,0){1.5}}
\put(2,0){\line(2,3){1}}
\put(6,0){\line(-2,3){1}}
\put(6,2){\line(2,-3){1}}
\put(6,2){\line(2,3){1}}
\put(8,2){\line(3,-2){1.5}}
\put(8,3){\line(1,0){0.5}}

\linethickness{0.1cm}
\put(2.67,0.5){\circle{0.04}}
\put(2.67,0.59){\makebox{0}}
\put(6.67,2.5){\circle{0.04}}
\put(6.47,2.1){\makebox{10}}

\linethickness{0.05cm}
\put(3.33,1.17){\circle{0.04}}
\put(3.34,1.25){\makebox{1}}
\put(3.66,1){\circle{0.04}}
\put(3.76,1.1){\makebox{2}}
\put(4.33,1){\circle{0.04}}
\put(4.13,1.1){\makebox{3}}
\put(4.67,1.17){\circle{0.04}}
\put(4.68,1.25){\makebox{4}}
\put(4.67,1.83){\circle{0.04}}
\put(4.74,1.70){\makebox{5}}
\put(5.33,1.5){\circle{0.04}}
\put(5.05,1.5){\makebox{6}}
\put(5.67,1.67){\circle{0.04}}
\put(5.72,1.7){\makebox{7}}
\put(6.33,1.17){\circle{0.04}}
\put(6.1,1.20){\makebox{8}}
\put(6.67,1.5){\circle{0.04}}
\put(6.5,1.6){\makebox{9}}
\end{picture}
\end{center}

\noindent
We thus can decompose $\chi_{U(\omega,\mathcal{A}_\omega)}-\chi_{U(\omega,\mathcal{A}_\omega).(x,\theta)}$ into the sum of the following differences
$$\chi_{U \big( \omega_i,\mathcal{A}_\omega.(x_i,\theta_i) \big)} - \chi_{U \big( \omega_{i+1},\mathcal{A}_\omega.(x_{i+1},\theta_{i+1}) \big)}$$
where $\omega_i$, $\omega_{i+1}$ are two tilings in $\Xi^2_\Delta$ obtained from $\omega$ by a direct isometry $(x_i,\theta_i)$, resp. $(x_{i+1},\theta_{i+1})$ and containing the origin on the punctuation of two simplices having a common edge.\\
We will thus focus on the difference
$$\chi_{U \big( \omega_1,\mathcal{A}_\omega.(x_1,\theta_1) \big)} - \chi_{U \big( \omega_2,\mathcal{A}_\omega.(x_2,\theta_2) \big)}$$
with $\omega_1=\omega.(x_1,\theta_1)$, $\omega_2=\omega.(x_2,\theta_2)$ and $(x_k,\theta_k) \in \mathcal{G} \Big( U(\omega,\mathcal{A}_\omega) \Big)$.\\

\noindent
We can then write $\omega_2 = \omega_1.(x_{12},\theta_{12})$ with $(x_{12},\theta_{12})=(x_1,\theta_1)^{-1}(x_2,\theta_2)$.\\
Thus : $\chi_{U \big( \omega_1,\mathcal{A}_\omega.(x_1,\theta_1) \big)} - \chi_{U \big( \omega_2,\mathcal{A}_\omega.(x_2,\theta_2) \big)}$ is equal to  
$$\chi_{U \big( \omega_1,\mathcal{A}_\omega.(x_1,\theta_1) \big)} - \chi_{U \big( \omega_1,\mathcal{A}_\omega.(x_1,\theta_1) \big).(x_{12},\theta_{12})}$$
and is therefore the differential of the characteristic function of $U(\omega_3,\mathcal{A}'_\omega)$ where $\omega_3$ is the tiling in $\Xi^1_\Delta$ with the origin on the punctuation of the common edge of the two above simplices and $\mathcal{A}'_\omega$ is the patch obtained from $\mathcal{A}_\omega$ by a direct isometry bringing the origin on this punctuation and taking the adequate orientation.\\
Thereby,  the generators of $H_{\Xi^2_\Delta}$ are images under $d^2_{PV}$ of elements of $C(\Xi^1_\Delta,\ZZ)$ and the reciprocal inclusion is proved.
\end{enumerate}

\begin{flushright}
$\square$
\end{flushright}

\noindent
By theorem \refg{isocechpimsner} and the fact that $\check{H}^3 ( \Omega ; \ZZ) \simeq \check{H}^2 ( \Omega/S^1 ; \ZZ)$ , we thus obtain :

\begin{cor}
The top integer \v{C}ech cohomology of the hull is isomorphic to the integer group of coinvariants of $\Xi_\Delta^2$ :
$$\check{H}^3 \big( \Omega ; \ZZ \big) \simeq C(\Xi_\Delta^2,\ZZ) / H_{\Xi_\Delta^2}.$$
\end{cor}

\noindent
We link now the coinvariants on the $\Delta$-transversal to the coinvariants of the canonical transversal which will be useful to prove the gap labelling in the next section.\\
The invariant ergodic probability measure $\mu$ on $\Omega$ induces a measure $\mu_0^t$ on each transversal $\Xi_0$ of the lamination which is given, locally, by : if $\big( V_i \times S_i \times C_i , h_i^{-1}\big)_i$ is a maximal atlas of the lamination and $B$ a borelian set in some $C_i$, we have 
$$\mu_0^t \big( B \big) = \dfrac{\mu \Big( h_i^{-1}(V_i \times S_i \times B) \Big)}{\lambda(V_i \times S_i)}$$
where $\lambda$ is a left and right Haar measure on $\RR^2 \times S^1$ (with $\lambda \big( [0;1]^2 \times S^1 \big) = 1$), $h_i^{-1}:V_i \times S_i \times C_i \longrightarrow U_i$, $V_i$ is an open subset of $\RR^2$, $S_i$ an open subset of $S^1$ and $C_i$ a clopen in $\Xi_0$.\\

\noindent
We then have a link between the $\Delta$-transversal and the canonical transversal : 

\begin{lem} \label{transversales}
Let $\mu_3^t$ be the induced measure on the transversal $\Xi \cup \Xi_\Delta^2$.\\
Denoting $\mu^t$ (resp. $\mu^t_2$) the restriction of $\mu^t_3$ to $\Xi$ (resp. $\Xi_\Delta^2$), we have : 
$$\mu^t_2\Big(C(\Xi_\Delta^2,\ZZ) \Big) \, = \, \mu^t \Big( C(\Xi,\ZZ) \Big).$$
\end{lem}

\textbf{Proof :}
		\begin{enumerate}
			\item[] In each prototile, there is $8$ punctuations of $2$-simplices which we number from $1$ to $16$ once for all (we are considering the $2$ prototiles (uncollared) of the pinwheel tiling and then we number the punctuations in the first prototile from $1$ to $8$ and those of the second prototile from $9$ to $16$). Thus, we can define the vectors $x_{0,i}$ joining the punctuation of the prototile taken for the definition of $\Xi$ (see \refg{espacedesprototuiles}) to the punctuation of the $i$-th simplex of this prototile.\\
Define the map $\Psi : C(\Xi_\Delta^2,\ZZ) \longrightarrow C(\Xi,\ZZ)$ by :
$$\Psi(f)(\omega_0)  =  \left\{ \begin{array}{l}
 {\displaystyle \sum_{i=1}^{8}} f(\omega_0-x_{0i})\\
 {\displaystyle \sum_{i=9}^{16}} f(\omega_0-x_{0i})\\
\end{array}\right.$$
depending on the tile type of the tile containing the origin in $\omega_0$.\\
This map send a function $f$ defined on $\Xi_\Delta^2$ on the function on $\Xi$ defined on a tiling $\omega_0$ containing the origin on the punctuation of a tile, by the sum of the values of $f$ on the tilings containing the origin on the punctuation of the cells constituting this tile.\\
This defines a group homomorphism which is trivially surjective.\\
Indeed, fix a simplex of each prototiles (for example, the number $1$ for the first prototile and $9$ for the second prototile), then we can define a section of $\Psi$,  $s:C(\Xi,\ZZ) \longrightarrow C(\Xi_\Delta^2,\ZZ)$, as follows 
$$s(f)(\omega_\Delta)  =  \left\{ \begin{array}{lc}
 f(\omega_\Delta + x_{0,1}) & \text{ if } \omega_\Delta \text{ has the origin in a tile  congruent to}\\
 & \text{ the first prototile and in a "1" simplex } \\
 f(\omega_\Delta + x_{0,9}) & \text{ if } \omega_\Delta \text{ has the origin in a tile  congruent to}\\
& \text{ the second prototile and in a "9" simplex}  \\
0 & \text{ else }
\end{array}\right.$$ 

\noindent
We then prove that this homomorphism preserves the measures.\\
Consider $\omega \in \Xi_\Delta^2$ (we suppose that the origin is on the punctuation of the $k$-th simplex) and $f=\chi_{U(\omega,\mathcal{A}_{\omega})}$ a generator of $C(\Xi_\Delta^2,\ZZ)$.\\
We then have $\Psi(f)=\chi_{U(\omega,\mathcal{A}_{\omega}) + x_{0,k}}$ and thus
$$\mu^t(\Psi(f))=\mu^t(U(\omega,\mathcal{A}_{\omega})+ x_{0,k})=\mu^t_3(U(\omega,\mathcal{A}_{\omega})+ x_{0,k})$$
since $\mu^t$ is the restriction of $\mu_3^t$ to $\Xi$.\\
Thereby
$\mu^t(U(\omega,\mathcal{A}_{\omega})+ x_{0,k})=\mu^t_3(U(\omega,\mathcal{A}_{\omega}))$ by invariance of $\mu^t_3$.\\
Finally, we obtain $\mu^t(U(\omega,\mathcal{A}_{\omega})+ x_{0,k}) = \mu^t_2(U(\omega,\mathcal{A}_{\omega}))$.\\
Thus, $\Psi$ preserves the measures and $\mu^t_2 \Big( C(\Xi_\Delta^2,\ZZ) \Big) \subset \mu^t \Big( C(\Xi,\ZZ) \Big)$.\\

\noindent
The reciprocal inclusion is obtained using the section $s$ : \\
if $\mu^t(f) \in \mu^t \Big( C(\Xi,\ZZ) \Big)$ then 
$$\mu^t(f) = \mu^t(\Psi \circ s (f) )= \mu^t_2(s(f)) \in \mu^t_2 \Big( C(\Xi_\Delta^2,\ZZ) \Big),$$
hence the reciprocal inclusion is proved together with the lemma.
\end{enumerate}

\begin{flushright}
$\square$
\end{flushright}

\section{Proof of the gap-labelling for pinwheel tilings and explicit computations}

\subsection{Proof of the gap-labelling for pinwheel tilings}

\noindent
To prove the gap-labelling for pinwheel tilings, i.e
$$[C_{\mu^t_2}] \Big( Ch_\tau \big( K_1(\Omega) \big) \Big) \subset \mu_2^t(C(\Xi_\Delta^2,\ZZ))= \mu^t(C(\Xi,\ZZ)), $$
we must now show that the inclucion of $C(\Xi,\ZZ)$ in $C(\Xi,\RR)$ is, at the level of cohomologies, the map 
$$r^*: \check{H}^3(\Omega;\ZZ) \longrightarrow H^3_\tau(\Omega)$$
induced by the inclusion of sheaves described in \citeg{MooSch} (see the diagram p.\pageref{diagcommu}).\\

\noindent
We first look at the lifting in $\check{H}^2(\Omega /S^1 ; \ZZ)$ of the generators of $C(\Xi_\Delta^2,\ZZ)/H_{\Xi_\Delta^2}$.\\
For this section, we consider the sequence $\{\B_l,f_l\}_{l \in \NN}$ defined in theorem \refg{limitprojective}.\\
$\left( \Xi_{l,\Delta}(\sigma_l) \right)_{\sigma_l \in \mathcal{S}^2_l,l \in \mathbb{N}}$ is then a base of neighborhoods of $\Xi_\Delta^2$ and the characteristic functions $\chi_{\sigma_l}$ thus span $C(\Xi^2_\Delta,\ZZ)$.\\
Fix one of these characteristic maps $\chi_{\sigma_l}$. This function is in fact, by definition, a function in $C(\Sigma_l^2,\ZZ)$ and so defines a class in the cohomology group $H^2_{PV}(\B_0;C(\Sigma_l,\ZZ))$ which is isomorphic to the simplicial cohomology group of $\B_l$, $H^2(\B_l;\ZZ)$. It's the class of the cochain which sends each characteristic map $\sigma$ on $\B_l$ on the integer $1$ if  $\sigma = \sigma_l$ and 0 else.\\
We would like to know  the image of this class under the isomorphism linking the simplicial cohomology to the \v{C}ech cohomology.\\
For this, in view of \citeg{Hoc} \textbf{5-5} and \textbf{8-2}, we see that this isomorphism is obtained by considering the coverings of $\B_l$ by open stars associated to the iterated barycenter decompositions of $\B_l$.\\
Denote $\mathcal{U}_n$ the covering of $\B_l$ by the open stars of the $n$-th barycenter decomposition of $\B_l$, noting that the open stars are the interior of the union of the $2$-simplices surrounding the vertices of the simplical structure of the $n$-th barycenter decomposition.\\
The second integer \v{C}ech cohomology of $\B_l$ is then the direct limit over the $\mathcal{U}_n$'s of $\check{H}^2(\mathcal{U}_n;\ZZ)$.\\
We need to shrink these open sets to end the proof. We don't take the open coverings by whole open stars, but we will consider large enough open subsets of these open stars (large enough in order to still cover the space). It suffices to take, for example, the open sets obtained from the open stars by rescaling them from the star vertex by a factor $\frac{5}{6}$. We call these open subsets, \textbf{pseudo-open stars}.\\
In this way, we make sure that the centre of a $2$-simplex is still in the intersection of the $3$ pseudo-open stars surrounding its vertices and thus the family $\mathcal{U}'_n$ of pseudo-open stars, is still an open cover of $\B_l$.\\
Moreover, the family $\mathcal{U}'_n$ is still a cofinal family of  coverings of $\B_l$ and thus 
$$\check{H}^2(\B_l;\ZZ) = \lim_{\longrightarrow} \check{H}(\mathcal{U}'_n;\ZZ).$$
The class of $\chi_{\sigma_l}$ in the integer group of coinvariants is thus sent on the class, in $\check{H}^2(\B_l;\ZZ)$, of the cochain $g_{\sigma_l}$ of $\check{H}(\mathcal{U}'_0;\ZZ)$ defined by taking the null section on the $2$-simplices $(U'_1,U'_2,U'_3)$ whose intersection is not contained in the simplex $\sigma_l$ and the constant section equal to $1$ on the $2$-simplex formed by the $3$ pseudo-open stars surrounding $\sigma_l$ (taking the same orientation as $\sigma_l$).\\
Let denote $U_1$, $U_2$ et $U_3$ the $3$ open stars surrounding $\sigma_l$, we then see that the previous cohomology class in $\check{H}^2(\B_l ; \ZZ)$ lifts in  $\check{H}^2_c(U;\ZZ):=\check{H}^2 \Big( U, \partial U;\ZZ \Big ) $ (where $U = U_1 \cap U_2 \cap U_3$)  on the cochain $h_{\sigma_l}$ equal to $1$ on the intersection of the $3$ pseudo-open stars surrounding $\sigma_l$ and vanishing elsewhere.\\
Now, if we look at the cohomology of $\Omega/S^1$, we see that the isomorphism in theorem \refg{isocechpimsner} sends the class of $\sigma_l$ on the class of $\p^*_l([g_{\sigma_l}])$ in $\check{H}^2(\Omega/S^1;\ZZ)$ which lifts on $\p^*_l([h_{\sigma_l}])$ in $\check{H}^2_c(U \times \Xi_{l,\Delta}(\sigma_l) ; \ZZ)$.\\
Indeed, we have the following commutative diagram:
$$\xymatrix{\check{H}_c^2(U;\ZZ) \ar[r] \ar[d]_{\p_l^*} & \check{H}^2(\B_l;\ZZ) \ar[d]^{\p_l^*}\\
\check{H}_c^2(U \times \Xi_{l,\Delta}(\sigma_l);\ZZ) \ar[r] & \check{H}^2(\Omega/S^1;\ZZ)
}$$
where the horizontal maps are induced by inclusion of the open sets $U$ and $U \times \Xi_{l,\Delta}(\sigma_l) = \p_l^{-1}(U)$ in $\B_l$ and $\Omega/S^1$ respectively.\\

\noindent
We thus fix an open set $U \times K$ where $U$ is an open subset of $\RR^2$ and $K$ a clopen in $\Xi_\Delta^2$.\\
We then have the following commutative diagram:\\
\\
$$\label{diagcommu}\scriptsize{\hspace*{-2.5cm} \xymatrix{ &  C(K,\ZZ) \ar@{^{(}->}[rrr] & & &  C(K,\RR) \ar[ddrr]^{\mu_2^t}\\
& \check{H}^2_c(\RR^2 \times K ; \ZZ) \ar[ddl] \ar[dd] \eq[u] \ar[rrr]^{r^*} \eq[dr] & & &  H^2_{\tau c}(\RR^2 \times K) \eq[u] \\ 
& & \check{H}^3_c(\RR^2 \times S^1 \times K ; \ZZ) \ar[rr]^{r^*} \ar[dd] \ar[ddr] & & H^3_{\tau c}(\RR^2 \times S^1 \times K) \ar[dd] \eq[u] \ar[rr]^{\hspace{1cm}[C_{\mu_2^t}]}  & & \RR \\
\check{H}^2(\Omega/S^1 ; \ZZ) \eq[r] & \check{H}_c^2\big( (\Omega \setminus F)/S^1 ; \ZZ\big) \eq[dr]\\
& &  \check{H}_c^3\big( \Omega \setminus F ; \ZZ\big) \eq[r] & \check{H}^3(\Omega;\ZZ) \ar[r]^{r^*} & H^3_\tau(\Omega) \ar[uurr]_{[C_{\mu_2^t}]}
}}$$

\noindent
The horizontal arrows, involving cohomologies, are the restriction map $r^*$ defined in \citeg{MooSch}.\\
$C(K,\ZZ) \hookrightarrow C(K,\RR)$ is the inclusion.\\
$\check{H}^2_c(\RR^2 \times K) \simeq C(K,\ZZ)$ and $\check{H}^2_{\tau c}(\RR^2 \times K) \simeq C(K,\RR)$ are Thom isomorphisms and so are natural with respect to sheaf maps.\\
The vertical arrows are induced by inclusion of an open subset in a space and are thus natural for sheaf maps.\\

\noindent 
Next, we note that the following diagram is formed by commutative diagrams:
\xymatrix{ \check{H}^2_c(\RR^2 \times K ; \ZZ) \eq[d] \ar[r]^{\otimes \RR} &  \check{H}^2_c(\RR^2 \times K ; \RR)  \ar[r]^{r^*} \eq[d] & H^2_{\tau c}(\RR^2 \times K) \eq[d] \\ 
\check{H}^3_c(\RR^2 \times S^1 \times K ; \ZZ) \ar[r]^{\otimes \RR} & \check{H}^3_c(\RR^2 \times S^1 \times K ; \RR)  \ar[r]^{r^*} & H^3_{\tau c}(\RR^2 \times S^1 \times K)
}\\

\noindent
The left vertical map is the Gysin isomorphism which become, by tensorising, the integration along the fibers $S^1$ in real cohomologies, this one being sent on the integration along the "longitudinal" fibers $S^1$ in longitudinal cohomologies.\\

\noindent
The two diagrams: \\
\\
$\scriptsize{\xymatrix{
& \check{H}^2_c(\RR^2 \times K ; \ZZ) \ar[dl] \ar[d] \\
\check{H}^2(\Omega/S^1 ; \ZZ) \eq[r] & \check{H}^2\big( (\Omega \setminus F)/S^1 ; \ZZ\big) 
}}$ \hspace*{0.5cm} and $\hspace{0.5cm} 
\scriptsize{\xymatrix{
 \check{H}^3_c(\RR^2 \times S^1 \times K ; \ZZ) \ar[d] \ar[dr] \\
\check{H}^3 \big( \Omega \setminus F ; \ZZ \big) \eq[r] & \check{H}^3\big( \Omega ; \ZZ \big) 
}}$\\
\\
are commutative because the open subsets $\RR^2 \times K \simeq U \times K$ ( resp. $U \times S^1 \times K$) is an open subset of $\Omega/S^1$ (resp. $\Omega$) included in $(\Omega \setminus F)/S^1$ ( resp. $\Omega \setminus F$), since tilings in  $F$ necessarily have the origin on an edge (see the patches surrounding the origin in tilings of $F$ below). 

\hspace*{-1.15cm}\begin{picture}(2.5,2)(0,0) 
\linethickness{0.02cm}
\put(1,1){\line(1,0){2}}
\put(1,1){\line(0,1){1}}
\put(1,2){\line(2,-1){2}}
\put(1,2){\line(1,0){2}}
\put(3,1){\line(0,1){1}}
\linethickness{0.075cm}
\put(2,1.5){\circle{0.05}}
\put(2,1.55){\makebox{O}}
\end{picture}
\begin{picture}(2.5,2)(0,0) 
\linethickness{0.02cm}
\put(1,1){\line(1,0){2}}
\put(1,1){\line(0,1){1}}
\put(1,1){\line(2,1){2}}
\put(1,2){\line(1,0){2}}
\put(3,1){\line(0,1){1}}
\linethickness{0.075cm}
\put(2,1.5){\circle{0.05}}
\put(1.75,1.55){\makebox{O}}
\end{picture}
\begin{picture}(3.25,3)(0,0.5) 
\linethickness{0.02cm}
\put(2,2){\line(1,0){2}}
\put(2,2){\line(0,1){1}}
\put(2,3){\line(2,-1){2}}
\put(1,2){\line(1,0){2}}
\put(1,2){\line(2,-1){2}}
\put(3,2){\line(0,-1){1}}
\linethickness{0.075cm}
\put(2.5,2){\circle{0.05}}
\put(2.5,2.1){\makebox{O}}
\end{picture}
\begin{picture}(3,3)(0,0.5) 
\linethickness{0.02cm}
\put(3,2){\line(-1,0){2}}
\put(3,2){\line(0,1){1}}
\put(3,3){\line(-2,-1){2}}
\put(4,2){\line(-1,0){2}}
\put(4,2){\line(-2,-1){2}}
\put(2,2){\line(0,-1){1}}
\linethickness{0.075cm}
\put(2.5,2){\circle{0.05}}
\put(2.5,2.1){\makebox{O}}
\end{picture}

\begin{picture}(5,5.5)(-2,-3) 
\linethickness{0.02cm}
\put(0,-2){\line(0,1){4}}
\put(-1,-2){\line(0,1){4}}
\put(1,-2){\line(0,1){4}}
\put(-1,2){\line(1,0){2}}
\put(-1,-2){\line(1,0){2}}
\put(-1,-2){\line(1,2){2}}
\put(-1,2){\line(1,-2){2}}
\put(-1,0){\line(1,0){2}}
\linethickness{0.075cm}
\put(0.2,0.1){\makebox{O}}
\end{picture}
\begin{picture}(5,5.5)(-3,-3) 
\linethickness{0.02cm}
\put(-2,0){\line(1,0){4}}
\put(-2,0){\line(0,-1){1}}
\put(-2,-1){\line(2,1){2}}
\put(2,0){\line(0,1){1}}
\put(2,1){\line(-2,-1){2}}
\put(-1,0){\line(0,1){2}}
\put(1,0){\line(0,-1){2}}
\put(-1,2){\line(1,-2){2}}
\put(1,-2){\line(-1,2){1}}
\put(0,0){\line(-3,-4){1.2}}
\put(0,0){\line(3,4){1.2}}
\put(-2,-1){\line(4,-3){.8}}
\put(2,1){\line(-4,3){.8}}
\put(-1,2){\line(4,-3){1.6}}
\put(1,-2){\line(-4,3){1.6}}
\linethickness{0.075cm}
\put(-0.15,-0.5){\makebox{O}}
\end{picture}

\noindent
Thus, the image in $H^3_\tau(\Omega)$ of a coinvariant generator $[\chi_{\sigma_l}]$ is sent under the Ruelle-Sullivan current in $\mu_2^t \big( C(K,\ZZ) \big) \subset \mu_2^t \big( C(\Xi_\Delta^2,\ZZ) \big) = \mu^t \big( C(\Xi,\ZZ) \big)$.\\

\noindent
Thanks to results obtained in \citeg{Hai}), we thus have
$$\tau^\mu_* \Big( K_0 \big( C(\Omega) \rtimes \RR^2 \rtimes S^1 \big) \Big) \subset \mu^t \Big( C(\Xi,\ZZ) \Big)$$

\noindent
The reciprocal inclusion is obtained by using the diagram in p.\pageref{diagcommu}.\\
Indeed, if you take a generator $\chi_{\sigma_l}$ in $C(\Xi_\Delta^2,\ZZ)$ it lifts to some $C(K,\ZZ)$ ($K$ is in fact $\Xi_{l,\Delta}(\sigma_l)$) which can be lifted using Gysin and Thom isomorphisms on a class in $\check{H}^3(\Omega;\ZZ)$ and since the Chern character is surjective, we have a lift $[u]$ of this class in $K_1 \big( C(\Omega) \big)$.\\
We then have 
$$\tau^\mu_*\big( \beta \circ \delta \circ \beta_{S^1} \circ \delta_{S^1} ([u]) \big) = [C_{\mu_2^t}] \big(r^*Ch([u]) \big) = \mu^t_2 \big( \chi_{\sigma_l} \big)$$
where $ \beta \circ \delta \circ \beta_{S^1} \circ \delta_{S^1} : K_1 \big( \Omega \big) \longrightarrow K_0 \big( C(\Omega) \rtimes \RR^2 \rtimes S^1 \big) $ is a map defined in \citeg{Hai} and is, in fact, the Kasparov product by the unbounded triple defined by the Dirac operator along the leaves of the foliated structure on $\Omega$.\\
Thus $\mu^t \big( C(\Xi,\ZZ) \big) = \mu^t_2 \big( C(\Xi_\Delta^2,\ZZ) \big) \subset \tau^\mu_* \Big( K_0 \big( C(\Omega) \rtimes \RR^2 \rtimes S^1 \big) \Big)$.\\

\noindent 
The gap-labelling is thus proved for pinwheel tilings:

\begin{theo}
If $\TT$ is a pinwheel tiling, $\Omega = \Omega(\TT)$ its hull provided with an invariant ergodic probability measure $\mu$ and $\Xi$ the canonical transversal provided with the induced measure $\mu^t$, we have :
$$\tau_*^\mu \Big( K_0 \big( C(\Omega) \rtimes \RR^2 \rtimes S^1 \big ) \Big) = \mu^t \big( C(\Xi,\ZZ) \big).$$
\end{theo}
\text{ }\\
\\

\subsection{Explicit computations}

\noindent
Now, we have $\tau_\mu^* \Big( K_0 \big( C(\Omega) \rtimes \RR^2 \rtimes S^1 \big) \Big) = \mu^t \big( C(\Xi,\ZZ) \big)$ and we want to compute $\mu^t \big( C(\Xi,\ZZ) \big)$ explicitly.\\

\noindent
For this, we will write $C(\Xi,\ZZ)$ as the direct limit of some system $(C(\Lambda_l,\ZZ),f^l )$ as we have done for $C(\Xi^2_\Delta,\ZZ)$ in \refg{limitdeltatransv}.\\
To define the group $C(\Lambda_l,\ZZ)$, we put a point in each tile of each $l$-prototile of $\B_l$.
\begin{center}
\begin{picture}(6,4)(0,0) 
\linethickness{0.05cm}
\put(1,1){\line(1,0){4}}
\put(1,1){\line(0,1){2}}
\put(1,3){\line(2,-1){4}}
\put(1,1){\line(1,2){0.8}}
\put(3,1){\line(-2,1){1.6}}
\put(3,1){\line(1,2){0.4}}
\put(3,1){\line(-3,4){1.2}}
\put(2,2){\circle{0.05}}
\put(1.2,2.4){\circle{0.05}}
\put(1.6,1,2){\circle{0.05}}
\put(2.8,1.6){\circle{0.05}}
\put(3.6,1.2){\circle{0.05}}
\end{picture}
\end{center}
Now, in the same way as we did for the $\Delta$-transversal, thanks to characteristic maps, if we take a tile $p$ in some supertile $\pi_j$ of $\B_l$, we denote $\Xi(p)$ the lifting of the punctuation of $p$ under $\p_l$, and $\chi_p$ its characteristic function in $\Xi$ (i.e $\chi_p([\omega])=1$ if and only if $\p_l([\omega]) = punct(p)$).\\
Thus, $\Xi(p)$ is formed by the tilings in $\Omega /S^1$ having the origin on the punctuation of a representative of $p$ which is itself contained in a representative of the supertile $\pi_j$.\\
We then define $C(\Lambda_l,\ZZ)$ as the subgroup of $C(\Xi,\ZZ)$ spanned by the characteristic functions $\chi_p$ where $p$ is a tile of some $l$-supertile of $\B_l$.\\
We also define $f^l : C(\Lambda_{l-1},\ZZ) \longrightarrow C(\Lambda_l,\ZZ)$ on the generators : let $p$ be a tile in some $(l-1)$-supertile $\pi_i$. We can uniquely decompose $\Xi(p)$ as the disjoint union $\Xi(p_j)$ where $p_j$ is a tile in some $l$-supertile with $f_l(punct(p_l)) = punct(p)$. $f^l(\chi_p)$ is then the sum of the $\chi_{p_j}$'s.\\

\noindent
We then see that ${\displaystyle C(\Xi,\ZZ) = \lim_{\longrightarrow} \big( C(\Lambda_l,\ZZ) , f^l \big)}$.\\
Indeed, $\Lambda_l$ is the set of all the punctuations in $\B_l$ and by theorem \refg{limitprojective}, we have $\Xi = \displaystyle{\lim_{\longleftarrow} (\Lambda_l,f_l)}$.\\

\begin{landscape}
\begin{figure}[ht] 
\begin{center}
\includegraphics[scale=0.25]{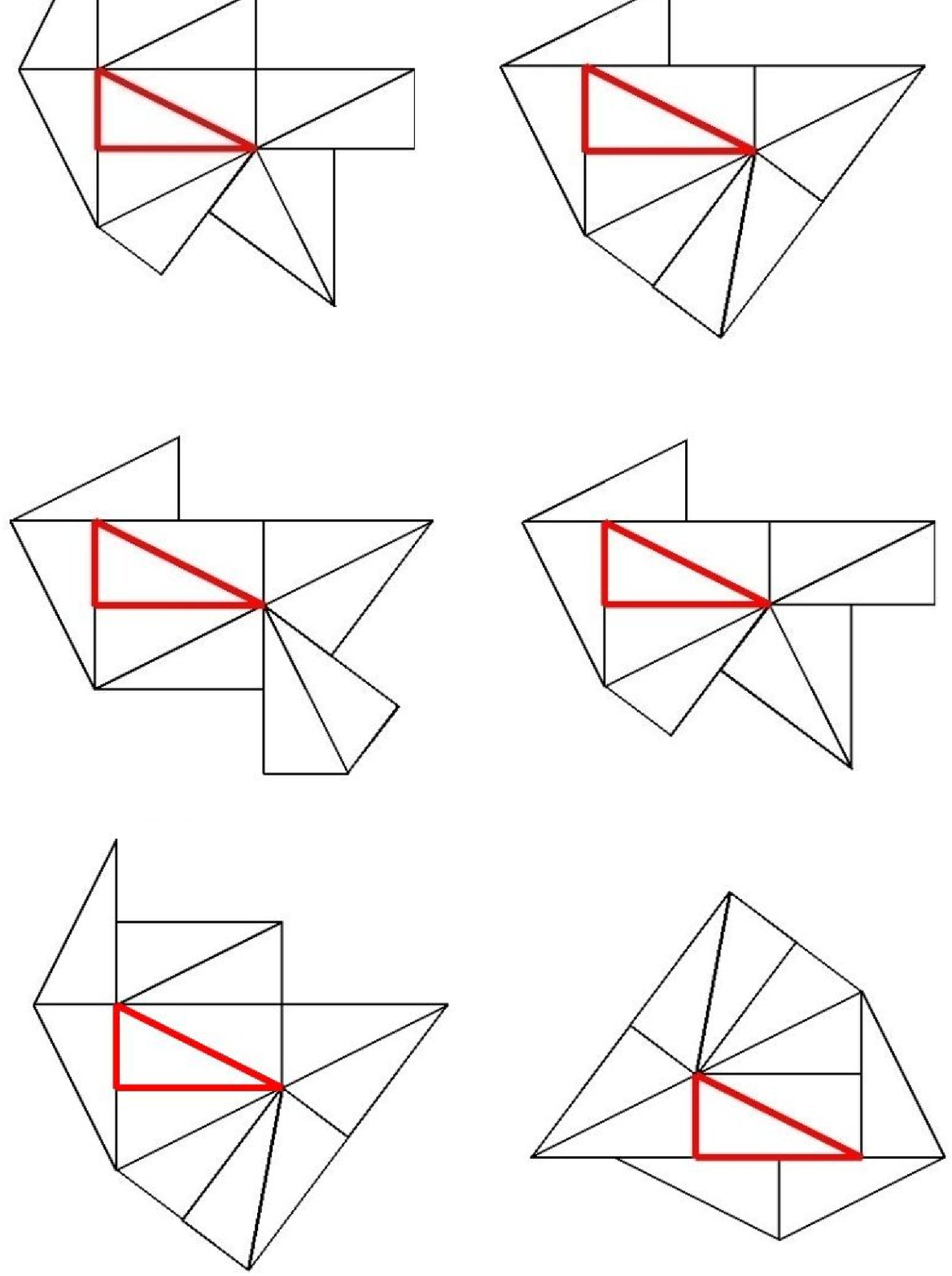}
\hspace{0.2cm}
\includegraphics[scale=0.25]{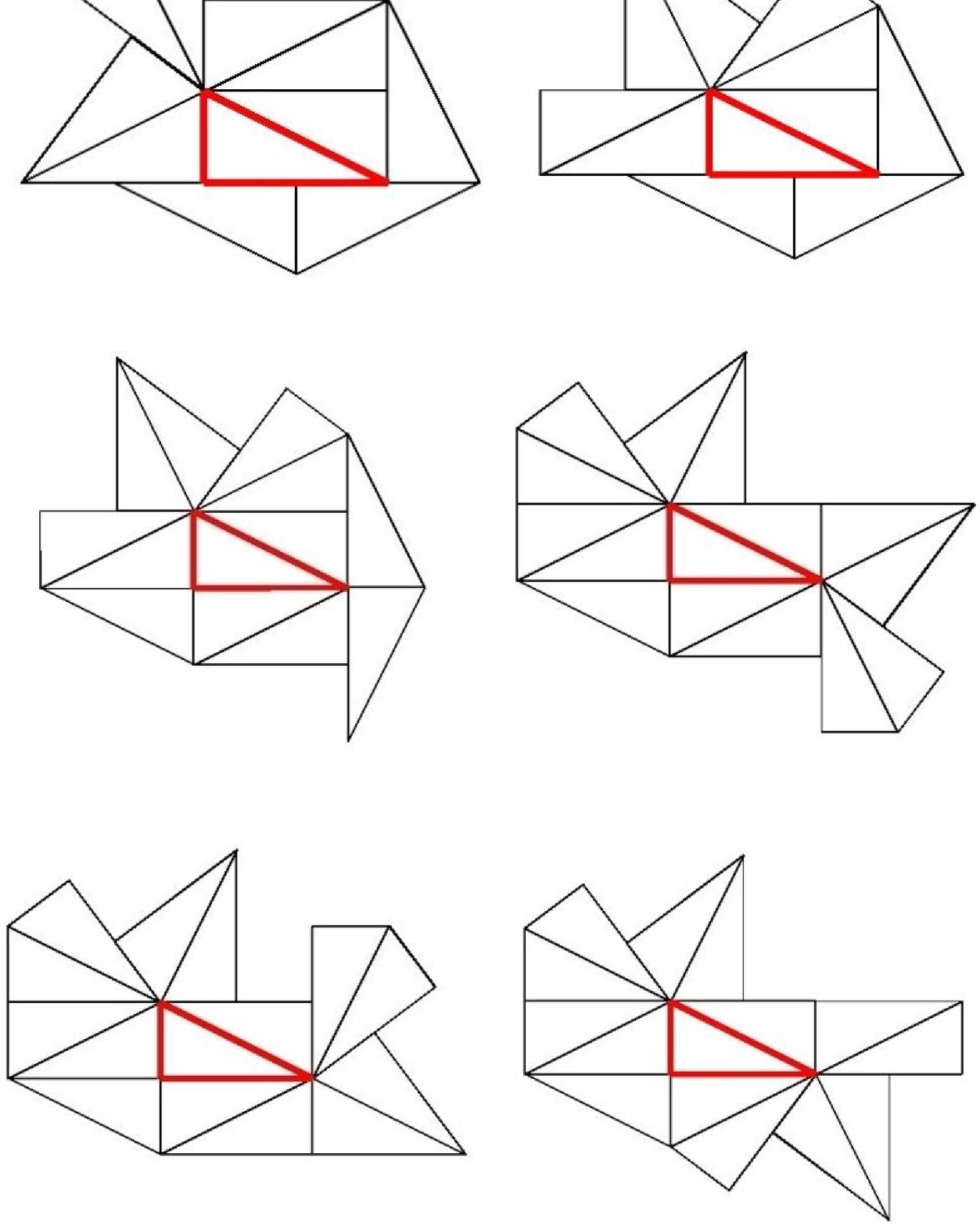}
\hspace{0.2cm}
\includegraphics[scale=0.25]{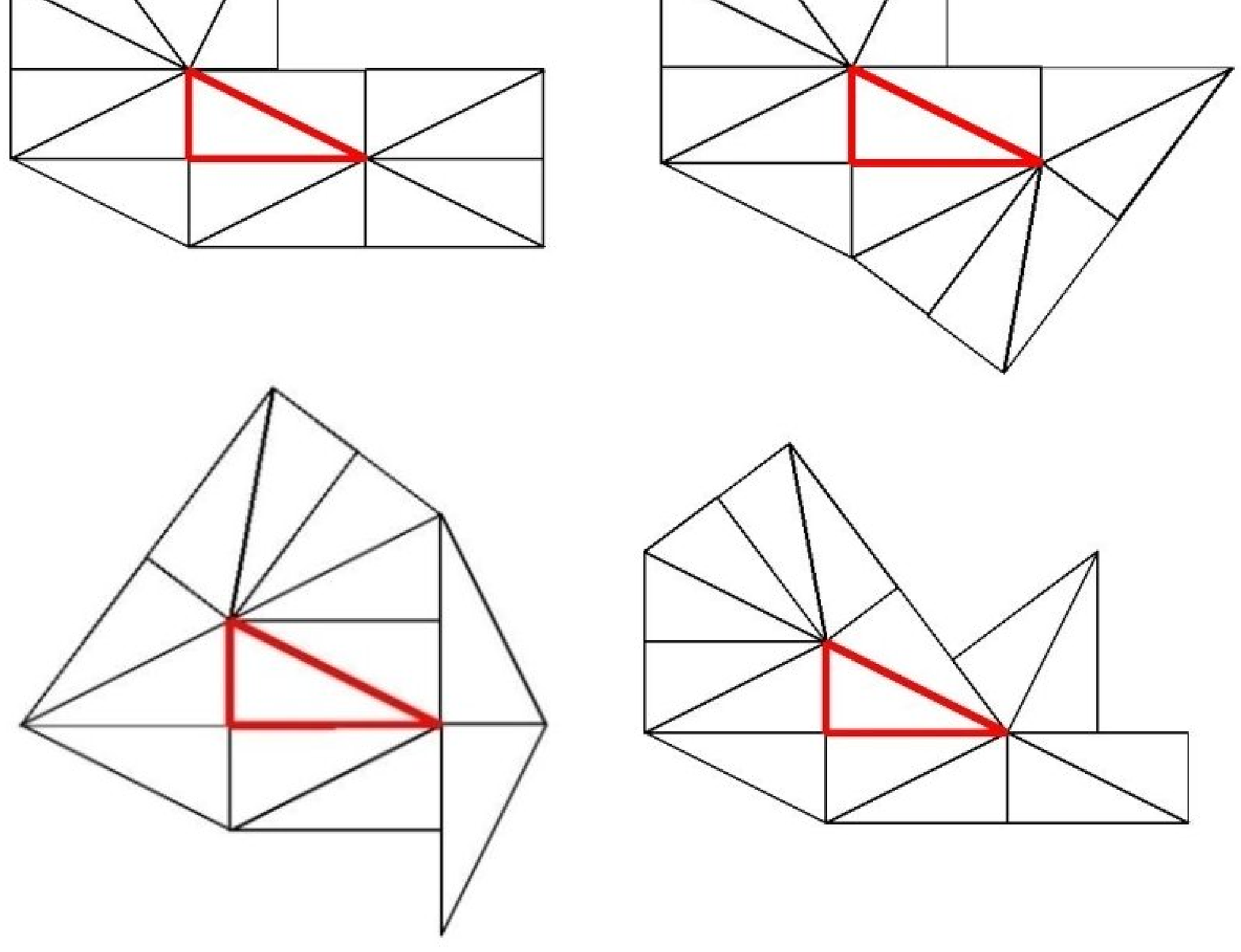}
\hspace{0.2cm}
\includegraphics[scale=0.25]{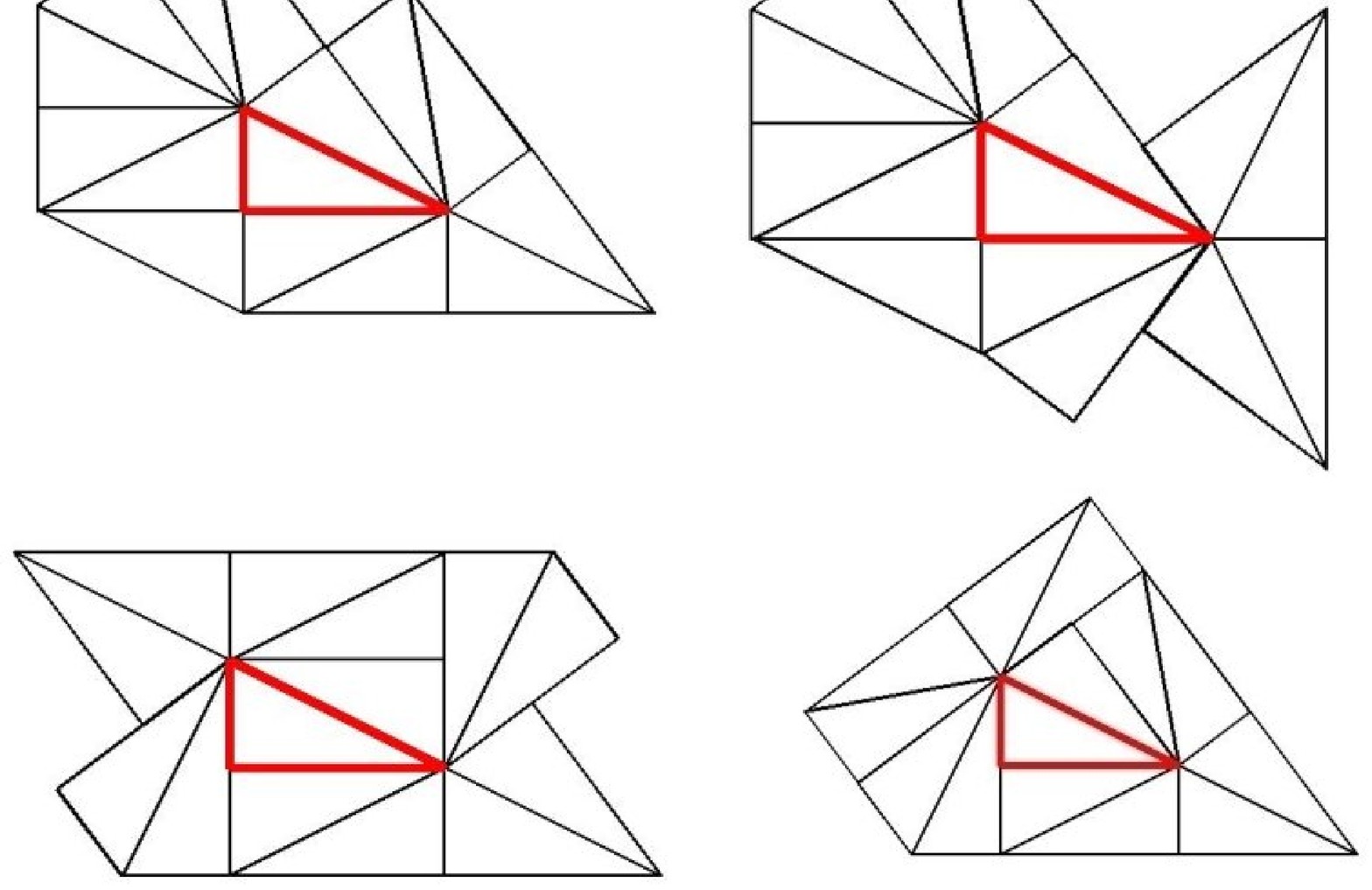}
\hspace{0.5cm}
\includegraphics[scale=0.25]{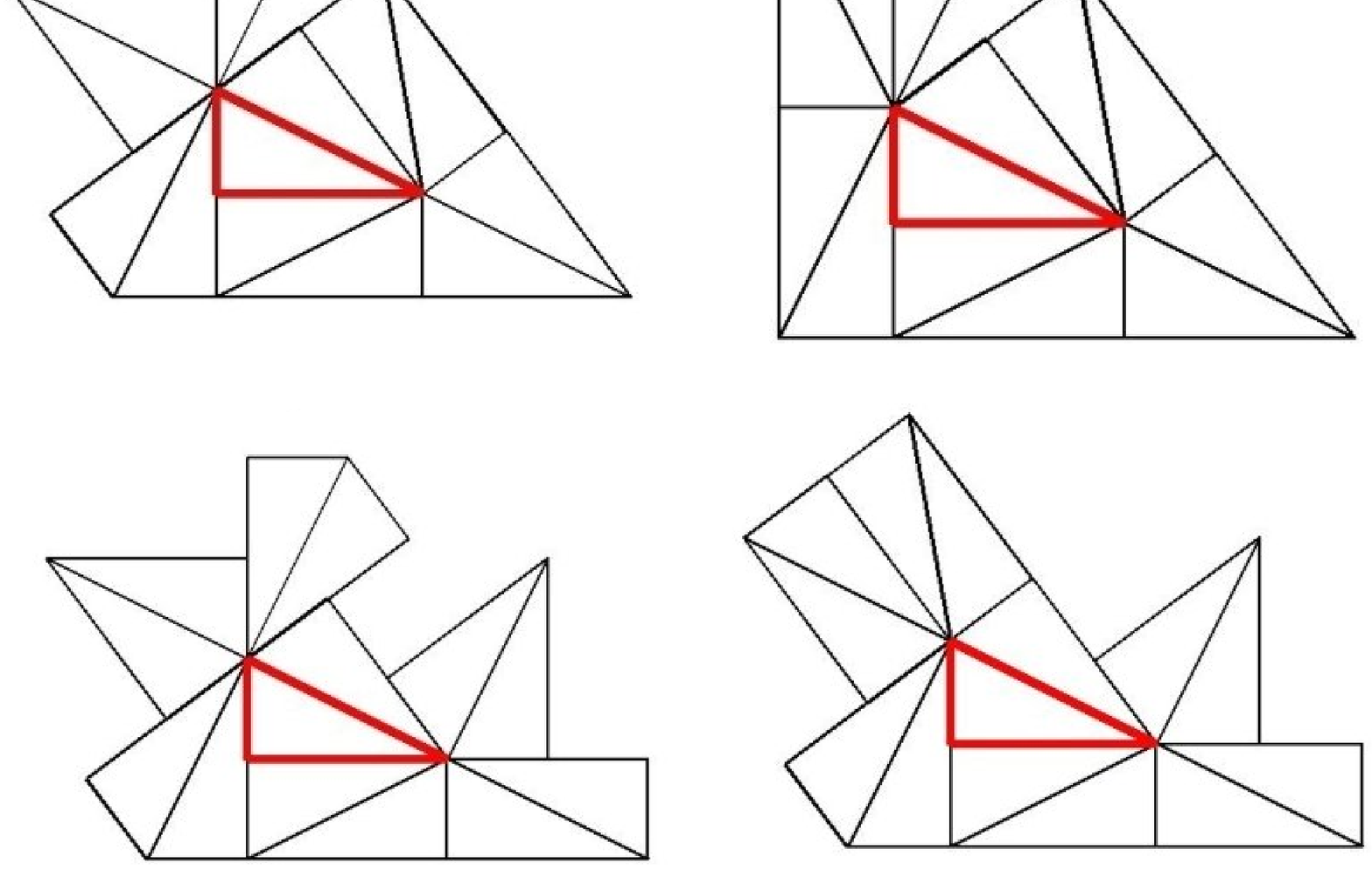}
\hspace{0.5cm}
\includegraphics[scale=0.25]{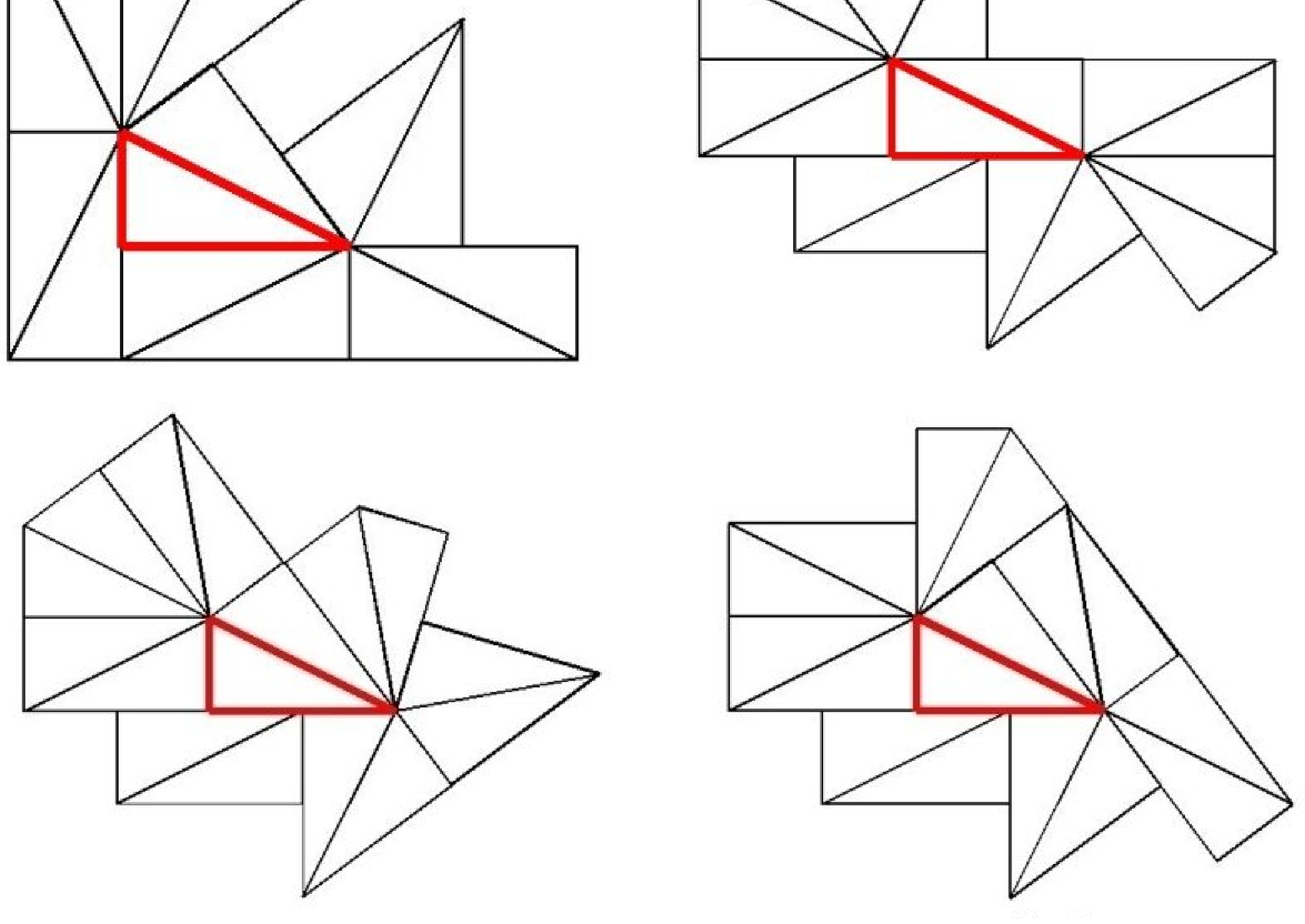}
\hspace{0.5cm}
\end{center}
\vspace{-3.5cm}\caption{28 collared prototiles. \label{prototuilescouronnees}} 
\end{figure}
\end{landscape}

\begin{landscape}
\begin{figure}[ht] 
\vspace{-1cm}
\begin{center}
\includegraphics[scale=0.25]{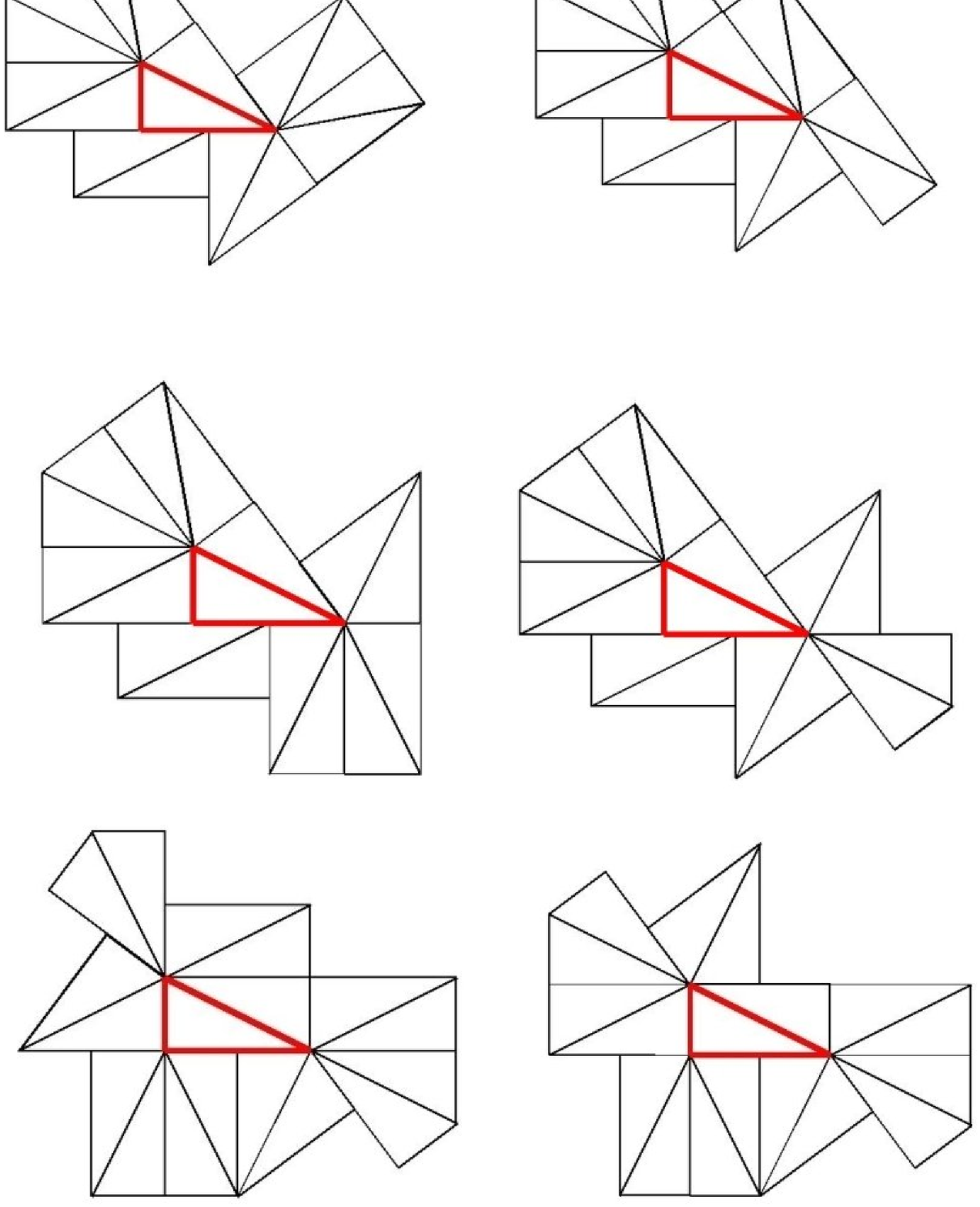}
\hspace{0.5cm}
\includegraphics[scale=0.25]{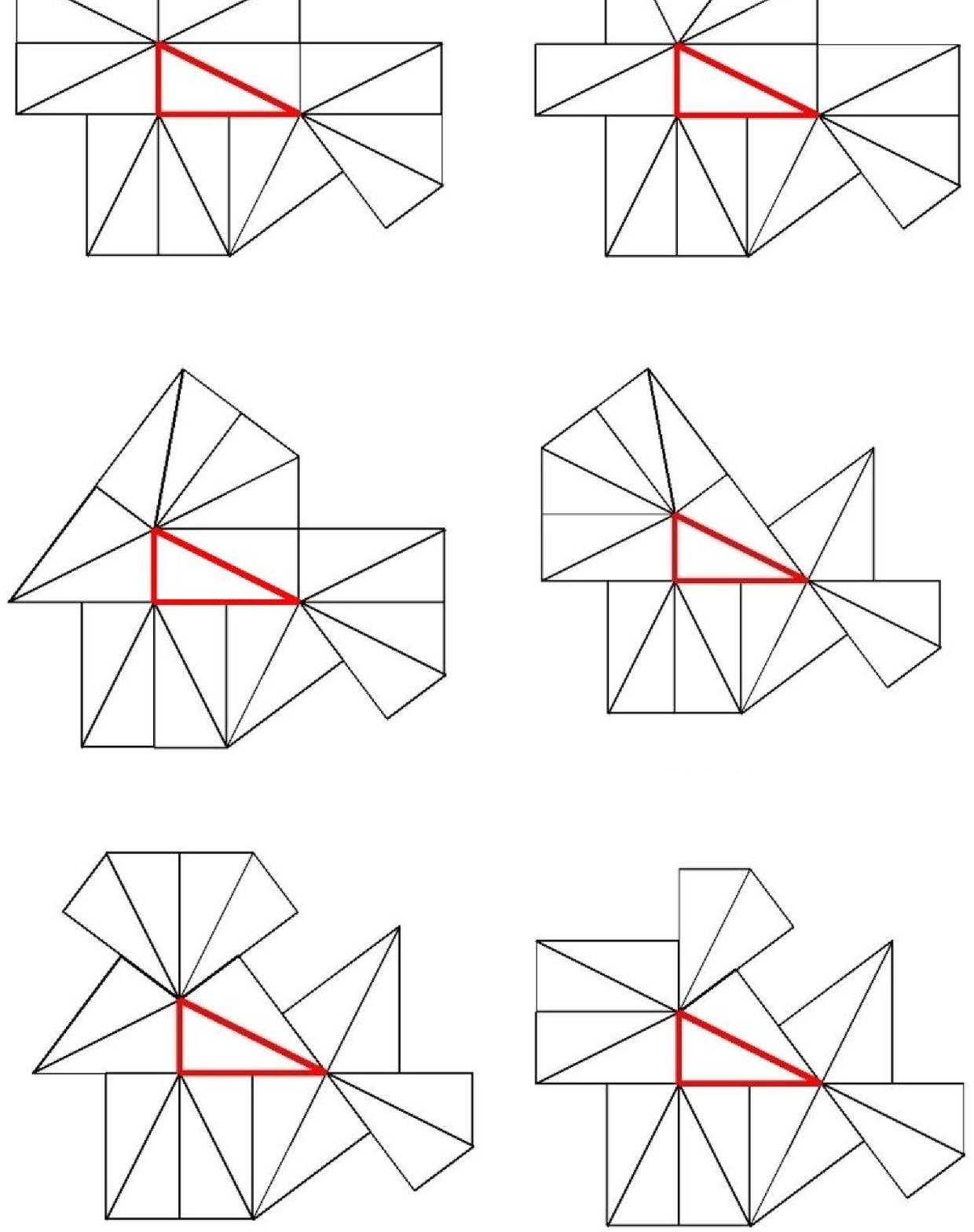}
\hspace{0.5cm}
\includegraphics[scale=0.25]{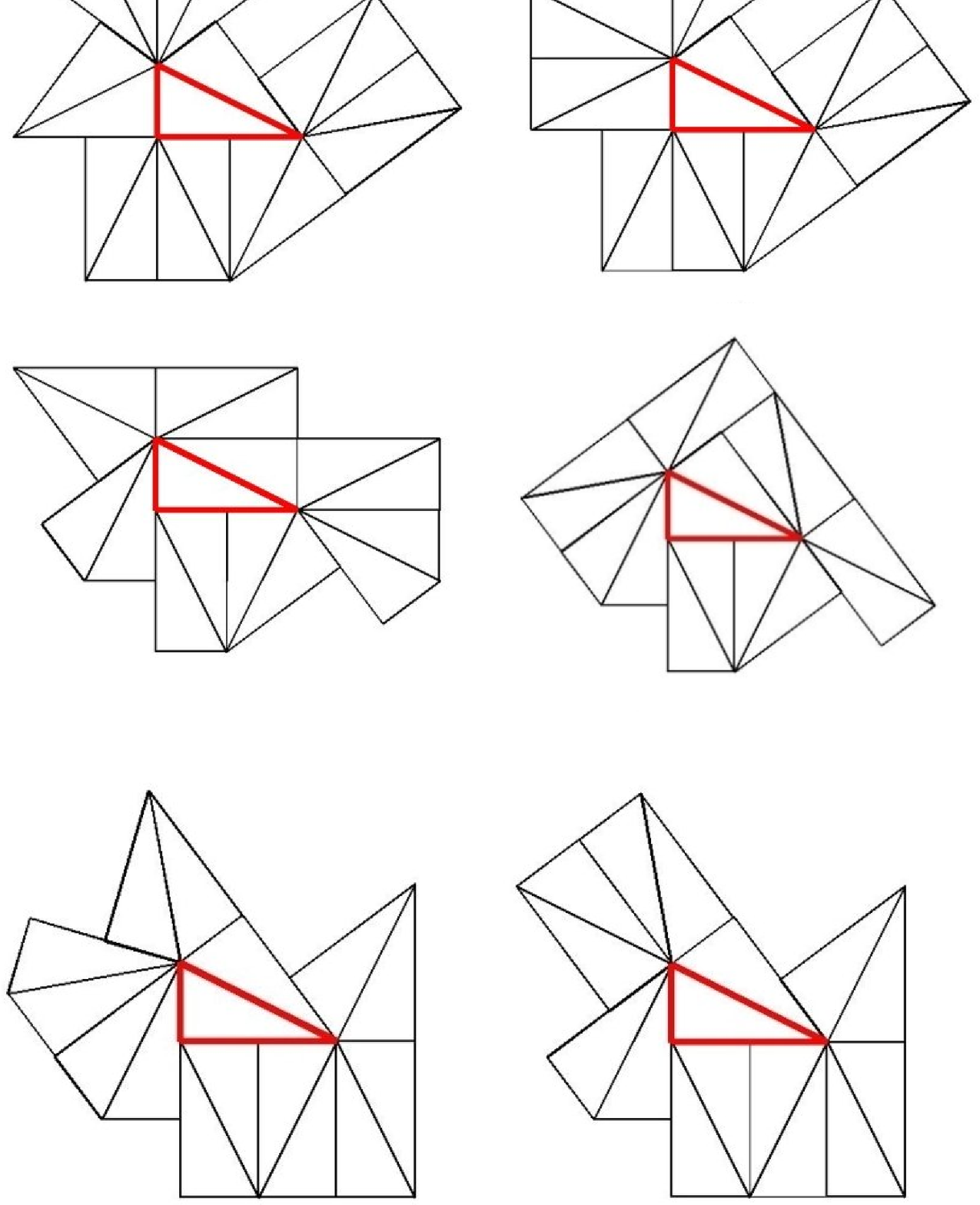}
\hspace*{0.5cm}
\includegraphics[scale=0.25]{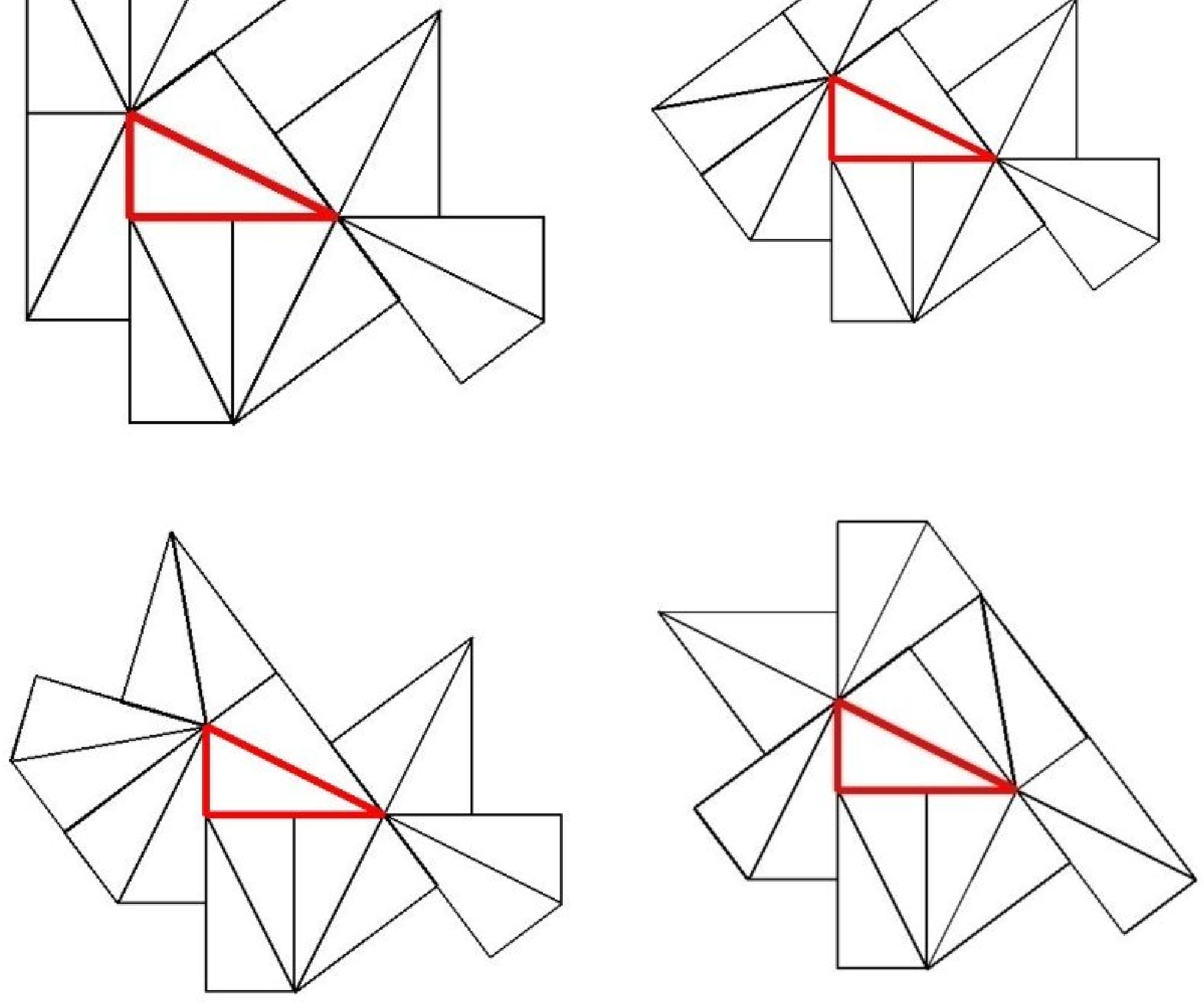}
\hspace*{0.5cm}
\includegraphics[scale=0.35]{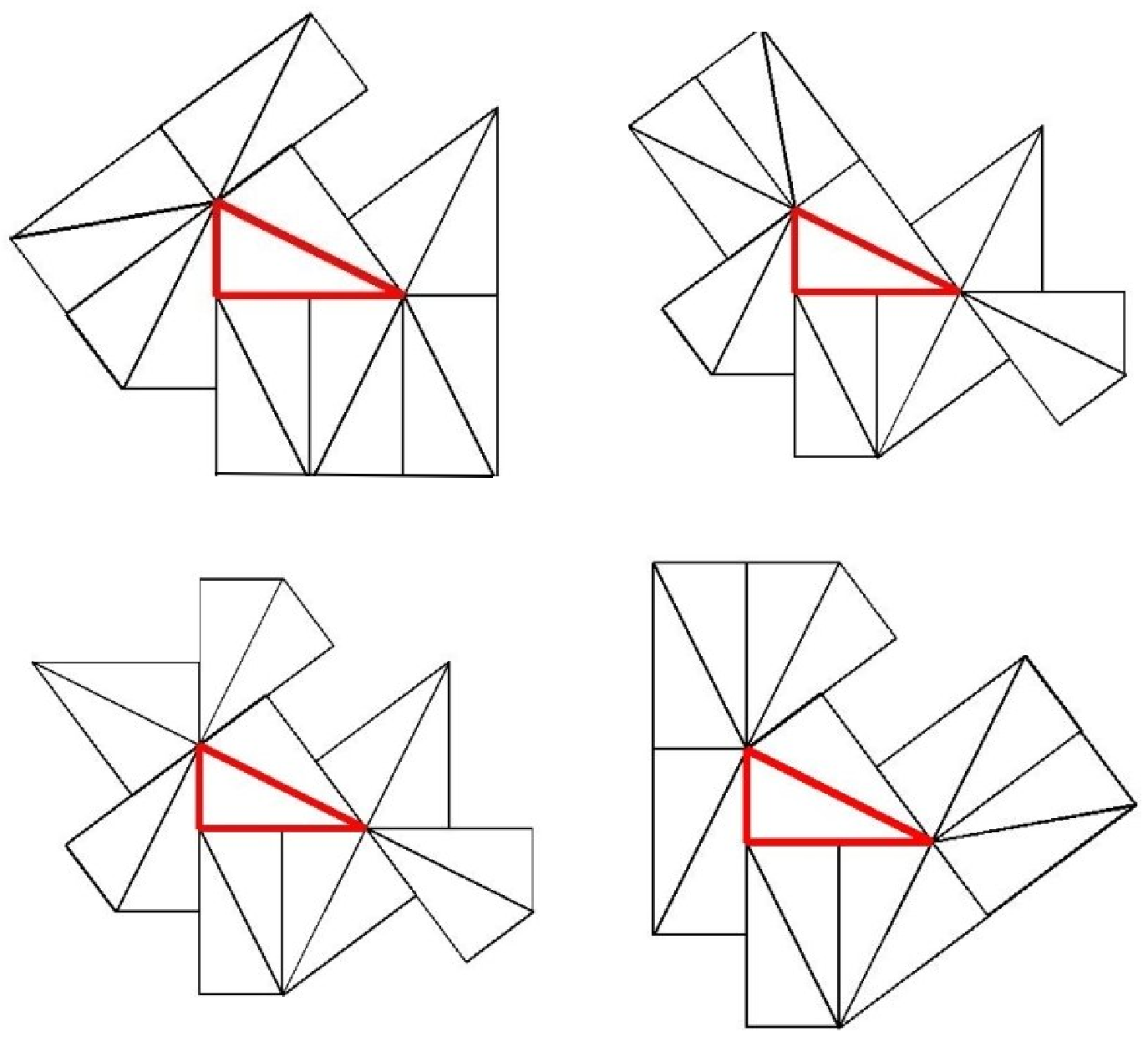}
\end{center}
\vspace{-2.5cm}\caption{26 other collared prototiles. \label{prototuilescouronnees2}}
\end{figure}
\end{landscape}

\noindent
For each $l \in \NN$, we can consider the subgroup $R_l$ of $C(\Lambda_l,\ZZ)$ spanned by differences $\chi_p - \chi_{p'}$ where $p$, $p'$ are two tiles in the same supertile of $\B_l$.\\
We then have, for each $l$, $C(\Lambda_l,\ZZ) / R_l \simeq \ZZ^{108}$, since there are exactly 108 collared supertiles in a pinwheel tiling (see \reffig{prototuilescouronnees} and \reffig{prototuilescouronnees2} on p.\pageref{prototuilescouronnees} and p.\pageref{prototuilescouronnees2}).\\

\noindent
Let $q_l:C(\Lambda_l,\ZZ) \longrightarrow C(\Lambda_l,\ZZ)/R_l$ be the quotient map.\\
$f^l$ factorizes through the quotient since, if $p$ and $p'$ are two tiles in the same $(l-1)$-supertile, $f^l(\chi_p - \chi_{p'}) = \sum \chi_{p_i} - \chi_{p'_i}$ where $p_i$ and $p_i'$ are two tiles contained in the same $l$-supertile of $\B_l$ and thus $f^l(\chi_p - \chi_{p'}) \in R_l$.\\
We then obtain a commutative diagram :
$$\xymatrix{C(\Lambda_{l-1},\ZZ) \ar[r]^{f^l} \ar[d]_{q_{l-1}} & C(\Lambda_l,\ZZ) \ar[d]^{q_l} \\
 C(\Lambda_{l-1},\ZZ)/R_{l-1} \ar[r]^{\tilde{f^l}} \eq[d] & C(\Lambda_l,\ZZ)/R_l \eq[d]\\
\ZZ^{108} \ar[r]^{A'} & \ZZ^{108}
}$$
where $A'$ is the transpose of the substitution matrix of the collared prototiles (i.e the matrix which have in position $(i,j)$ the number of representatives of the collared prototile of type $i$ in the substitution of the collared prototile of type $j$).\\

\noindent
The system $\Big( C(\Lambda_l,\ZZ)/R_l, \tilde{f^l} \Big)$ is a direct system and we can thus consider the direct limit $CR:={\displaystyle \lim_{\longrightarrow} \Big( C(\Lambda_l,\ZZ)/R_l, \tilde{f^l} \Big) }$.\\
Define also $\mathcal{C}:=C(\Xi,\ZZ) / \mathcal{R}$ where $\mathcal{R}$ is the subgroup of $C(\Xi,\ZZ)$ spanned by the $R_l$'s and denote $q:C(\Xi,\ZZ) \longrightarrow \mathcal{C}$ the quotient map.\\

\begin{lem}
The group $CR$ is isomorphic to $\mathcal{C}$.
\end{lem}

\noindent
\textbf{Proof :}
		\begin{enumerate}
			\item[] Let $l \in \NN$, we thus have a homomorphism $\psi_l:C(\Lambda_l,\ZZ)/R_l \longrightarrow \mathcal{C}$ defined by $\psi_l(q_l(f)):=q(f)$ ($R_l \subset \mathcal{R}$).\\
We have also the following commutative diagram :
$$\xymatrix{C(\Lambda_{l-1},\ZZ)/R_{l-1} \ar[rr]^{\tilde{f^l}} \ar[rd]_{\psi_{l-1}} & & C(\Lambda_l,\ZZ)/R_l \ar[dl]^{\psi_l} \\
 & \mathcal{C} &
}$$ 
since the next diagram is formed by commutative diagrams:\\
$\hspace*{0cm} \scriptsize{\xymatrix{C(\Lambda_{l-1},\ZZ)/R_{l-1} \ar[rrrr]^{\tilde{f^l}} \ar@/_2pc/[rrdddd]& & & & C(\Lambda_l,\ZZ)/R_l \ar@/^2pc/[lldddd]\\
& C(\Lambda_{l-1},\ZZ) \ar[ul]^{q_{l-1}} \ar[rr]^{f^l} \ar[ddr]_{id} & & C(\Lambda_l,\ZZ) \ar[ur]_{q_l} \ar[ddl]^{id} \\
\\
& & C(\Xi,\ZZ) \ar[d]^q & & \\
& & \mathcal{C} & & 
}}$\\

\noindent
Denoting $j_l:C(\Lambda_l,\ZZ)/R_l \longrightarrow CR$ the canonical homomorphisms sending an element of $C(\Lambda_l,\ZZ)/R_l$ on its class in the limit $CR$, there exists a unique homomorphism $j:CR \longrightarrow \mathcal{C}$ such that $j \circ j_l = \psi_l$ (by definition of direct limit).\\
We then prove that $j$ is an isomorphism.\\

\noindent
\underline{Surjectivity} : consider $q(f) \in \mathcal{C}$, then, since $C(\Xi,\ZZ)$ is the direct limit of the $C(\Lambda_l,\ZZ)$'s, $q(f)=q(f_i)$ with $f_i \in C(\Lambda_{l_i})$.\\
We then have $q(f)=q(f_i)=\psi_{l_i} \big( q_{l_i}(f_i) \big) = j \circ j_{l_i} \circ q_{l_i} (f_i)$ and thus $j$ is surjective.\\

\noindent
Now, we show that $j$ is also \underline{injective}.\\
Suppose that $j(f)=0$ ($f \in CR$), then, $f$ can be written $j_n(q_n(g))$ with $g \in C(\Lambda_n,\ZZ)$ (by definition of direct limit) and thus $j \circ j_n \circ q_n (g)=0$ i.e $\psi_n \circ q_n(g)=0$.\\
Hence, $g \in \mathcal{R}$ and we can find $k_1 \leqslant \ldots \leqslant k_r$, $c_1,\ldots, c_r \in \ZZ$ and $f_1, \ldots, f_r$ with $f_i \in R_{k_i}$ such that $g=c_1 f_1 + \ldots + c_r f_r$.\\
By definition, if $h \in C(\Lambda_j,\ZZ)$, then $f^{j+1}(h)=h$ in $C(\Xi,\ZZ)$ ( $h$ is written as a linear combination of characteristic functions of patches formed by supertiles of level $j$ and $f^{j+1}(h)$ is just another way to write this sum obtained by decomposing these patches in patches of $(j+1)$-supertiles that contain them).\\
Hence, $f_i=f_{k_rk_i}(f_i)$ for $i=1,\ldots,r-1$, with $f_{km}:=f^{k} \circ \ldots \circ f^{m+1}$ and
$$f_{k_r n}(g) = g = c_1 f_{k_rk_1}(f_1) + \ldots + c_{r-1} f_{k_r k_{r-1}}(f_{r-1}) + c_r f_r$$
in $C(\Xi,\ZZ)$ and $f_{k_r n}(g)=g$ is then in $R_{k_r}$.\\
Thus $\tilde{f_{k_r n}}(q_n(g))=q_{k_r}(f_{k_r n}(g))=0$ with $\tilde{f_{k_r n}}:=\tilde{f^{k_r}} \circ \ldots \circ \tilde{f^{n+1}}$.\\
We then conclude that  $j_{k_r} \circ \tilde{f_{k_r n}} \circ q_n(g)=0$ and by definition of direct limit, $0=j_{k_r} \circ \tilde{f_{k_r n}} \circ q_n(g) = j_n \circ q_n(g) = f$ and $j$ is thus injective.
		\end{enumerate}
\begin{flushright}
$\square$
\end{flushright}

\noindent
The previous lemma is not surprising since we can remark that $\mathcal{R}=\cup R_l$ and $C(\Xi,\ZZ) = \cup C(\Lambda_l,\ZZ)$.\\

\noindent
Moreover, $\mu^t(\chi_p-\chi_{p'} ) =0$ if $p$, $p'$ are two tiles in the same supertile since, if $\chi_p=\chi_{U(\omega,ST_\omega)}$ where $[\omega] \in \Xi(p)$, $ST_\omega$ being the supertile surrounding $p$ and such that the representative of $p$ containing the origin in $\omega$ has the "good" orientation, then  $\chi_{p'}=\chi_{U(\omega,ST_\omega).(x,\theta)}$ with $(x,\theta)$ the direct isometry sending $\omega$ on the tiling of $\Xi(p')$ having a representative of $p'$ surrounding the origin and in the "good" orientation.\\
Hence, the measure factorizes through the quotient and $\mu^t(\mathcal{C}) = \mu^t(C(\Xi,\ZZ))$.\\
We must then compute $\mu^t(\mathcal{C})$ to end the computation of the gap-labelling of the pinwheel tiling.\\
\\
Since the following diagram 
$$\xymatrix{ \cdots \ar[r] & C(\Lambda_{l-1},\ZZ) / R_{l-1} \ar[r]^{\tilde{f^l}} \eq[d] & C(\Lambda_l,\ZZ) / R_l \ar[r] \eq[d]& \cdots \\
\cdots \ar[r] & \ZZ^{108} \ar[r]_{A'} & \ZZ^{108} \ar[r] & \cdots
}$$
is commutative, we have $\mathcal{C} \simeq {\displaystyle \lim_{\longrightarrow}(\ZZ^{108},A')}$.\\

\noindent
We then use results obtained by Effros in \citeg{Eff} since we have a stationary system:
\begin{defi} (see \citeg{Eff})
A direct system
$$\xymatrix{\ZZ^r \ar[r]^\phi &  \ZZ^r \ar[r]^\phi  & \cdots
}$$
with constant spaces $\ZZ^r$ and with constant map $\phi$ (constant mean that it's always the same map) is called \textbf{stationary}.\\
A stationary system is \textbf{simple} if, for some $n$, $\phi^n$ is strictly positive (i.e all its coefficients are stricly positive).\\

\noindent
An \textbf{ordered group} $G$ is an abelian group $G$ together with a subset $P$, called the \textbf{positive cone} and denoted $G^+$, such that:
\begin{enumerate}
	\item $P + P \subset P$,
	\item $P - P = G$,
	\item $P \cap (-P) = \{0\}$,
	\item if $a \in G$ and $n a \in P$ for some $n \in \NN$, then $a \in P$.
\end{enumerate}
We shall write $a \leqslant b$ (resp. $a < b$) if $b-a \in G^+$ (resp. $G^+ \setminus \{0\}$).\\ 
We shall say that $u \in G^+$ is an \textbf{order unit} on $G$ if 
$$\{a \in G : 0 \leqslant a \leqslant n.u \text{ for some } n\in \NN \}=G^+.$$
A \textbf{state} (depending on $u$) on $G$ is a homomorphism $p:G \rightarrow \RR$ such that $p(G^+) \geqslant 0$ and $p(u)=1$.\\
We denote $S_u(G)$ the set of such states.
\end{defi}
For example, if $G$ is the limit of the system $\xymatrix{\ZZ^r \ar[r]^\phi &  \ZZ^r \ar[r]^\phi &  \cdots
}$ (with $\phi$ having all its coefficients non-negative), we can define an ordering on $G$ taking $G^+$ to be the union of all the images in $G$ of the $(\ZZ^+)^r$'s.\\

\noindent
With these definitions, we have the following theorem :

\begin{theo} (\citeg{Eff}) \text{}\\
Let $u_1$ be $(1,\ldots,1)$ in the first copy of $\ZZ^r$ in the simple stationary system of the definition then $u=\phi_\infty(u_1)$, its image in the limit $G$, is an order unit for $G$ and $S_u(G)$ contains only one point.
\end{theo}

\noindent
Our inverse system $(\ZZ^{108},A')$ is in fact stationary and simple (since $A$ is the substitution matrix of our tiling which is a primitive substitution ($A^6$ has all its coefficients strictly positive)) and thus there exists only one state $p$ on the inverse limit $G={\displaystyle \lim_{\longrightarrow}(\ZZ^{108},A')}$.\\
Moreover, we have an explicit formula for this state : if $\lambda$ is the Perron-Frobenius eigenvalue of $A'$ and $\alpha$ the unique eigenvector of $A$  associated to $\lambda$ such that $\sum \alpha_i = 1$, then
$$p([a,n]) = \frac{1}{\lambda^{n-1}}\sum \alpha_i a_i.$$

\noindent
In the case of pinwheel tilings, we have $\lambda=5$, $u_1=(1, \ldots,1) \in \ZZ^{108}$ and there is only one state on $G$ defined by :
$$p\Big( \big[ (k_1,\ldots,k_{108}) ; n \big] \Big) = \frac{1}{33000} \frac{1}{5^{n-1}} \sum_{i=1}^{108} k_i \alpha'_i,$$
with $\alpha = \dfrac{1}{33000} \alpha'$ and  
$$\begin{array}{rcl} 
\alpha' & = & 
(765,1185,360,255,735,1185,360,255,765,90,90,255,250,80,\\
 & & 735,255,80,400,360,660,600,660,300,360,360,255,400,360,\\
 & & 90,255,90,350,90,255,250,255,80,80,163,18,237,72,90,360,\\
 & & 237,72,204,204,163,300,51,18,50,51,765,1185,360,255,735,\\
 & & 1185,360,255,765,90,90,255,250,80,735,255,80,400,360,660,\\
 & & 600,660,300,360,360,255,400,360,90,255,90,350,90,255,250,\\
 & & 255,80,80,163,18,237,72,90,360,237,72,204,204,163,300,51,\\
 & & 18,50,51)
\end{array}$$

\noindent
Moreover, since each $C(\Lambda_l,\ZZ) / R_l$ is naturally ordered by the ordering of $\ZZ^{108}$, $\mathcal{C}$ have an ordering compatible with its isomorphism with $G$ and the unit order for this ordering is the class of the constant function equal to $1$.\\
Hence, there is a unique state on $\mathcal{C}$ such that, if $p_i$ are tiles contained in $108$ different $n$-supertiles of $\B_n$, then :
\begin{eqnarray} \label{pp}
p\Big( \sum k_i q_n(\chi_{p_i}) \Big) = \dfrac{1}{33000} \dfrac{1}{5^n} \sum_{i=1}^{108} k_i \alpha'_i.
\end{eqnarray}
But $\dfrac{\mu^t}{\mu^t(\Xi)}$ is in fact a state on $\mathcal{C}$ thus $\mu^t = \mu^t(\Xi)p$ and 
$$\mu^t(C(\Xi,\ZZ)) = \frac{\mu^t(\Xi)}{264} \ZZ \Big[\frac{1}{5} \Big].$$

\noindent
We end this article by two points.\\
\\
The first one is that, in fact, $\mu^t$ is a probability measure on $\Xi$ :

\begin{lem}
$\mu^t(\Xi) = 1$.
\end{lem}

\noindent
\textbf{Proof :}
		\begin{enumerate}
			\item[] This result is obtained using the formula in \citeg{MooSch} p.90 : if $f:\Omega \rightarrow \Xi$ is a Borel function with $f(x)$ in the same leaf as $x$, then for each point $\omega \in \Xi$ we can define $\rho_\omega$ on $f^{-1}(\omega)$ as the restriction of the Haar measure on the leaf $l(\omega)$ of $\omega$ to $f^{-1}(\omega) \subset l(\omega)$. We then have :
$$\mu(E) = \int_\Xi \left(  \int_{f^{-1}(\omega)} \chi_E(x) d\rho_\omega(x) \right) d\mu^t(\omega),$$
for any Borel set $E$ in $\Omega$.\\
\\
We then just need to produce such a Borel function.\\
For this, consider $\hat{p_i^c}$ ($i=1,\ldots,108$) the $108$ collared prototiles of the pinwheel tiling and let $V_i^k$ be the open subset of $\Omega$  homeomorphic to $\Xi(p_i^c) \times \overset{\quad\!\!\!\!\!\!\!\! \circ}{p_i^c} \times S_k$ where $ \overset{\quad\!\!\!\!\!\!\!\! \circ}{p_i^c}$ is the interior of the representative of $\hat{p_i^c}$ with the origin on its punctuation (and with the "good" orientation) and \\$S_k=\left] \frac{2k\pi}{3}-\frac{\pi}{3} , \frac{2k\pi}{3}+\frac{\pi}{3} \right[$ for $k=0,1,2$.\\
Set ${\displaystyle U_i = \bigcup_{k=0}^2 V_i^k}$ for $i=1, \ldots, 108$.\\
The union of the closure of these open sets is a covering of $\Omega$ :
$$\Omega = \bigcup_{i=1}^{108} \overline{U_i} = \bigcup_{i=1}^{108} \bigcup_{k=0}^2\overline{V_i^k}.$$

\noindent
Setting 
$$F^0_1=\overline{V_1^0},\, F^1_1=\overline{V_1^1} \setminus \overline{V_1^0}, \, F^2_1= {\displaystyle \overline{V_1^2} \setminus \left( \overline{V_1^0} \cup \overline{V_1^1} \right)},\, F^0_2=\overline{V^0_2} \setminus \overline{U_1},$$ 
$$F^1_2={\displaystyle \overline{V^1_2} \setminus \left( \overline{V_2^0} \cup \overline{U_1} \right) },\, \ldots,\, {\displaystyle F^2_{108}=\overline{V^2_{108}} \setminus \left( \bigcup_{l=1}^{107} \overline{U_l} \bigcup \overline{V^0_{108}} \bigcup \overline{V^1_{108}} \right) },$$ we have a covering of $\Omega$ by disjoint Borel sets $F_i$.\\
We then define $f:\Omega \rightarrow \Xi$ by $f(\omega)=\omega_i$ if $\omega \in F_i^k$, $\omega=\omega_i.(x_i,\theta_i)$ with $\omega_i \in \Xi(p_i^c)$ and $(x_i,\theta_i) \in  p_i^c \times \overline{S_k}$.\\
This is a Borel function from $\Omega$ to the transversal $\Xi$ and $f(\omega)$ is on the same leaf as $\omega$.\\
\\
We then have proved that $\mu$ can be reconstructed from $f$ and $\mu^t$ : 
\begin{eqnarray} \label{mu}
\mu(E) = \int_\Xi \left(  \sum_{k=0}^2 \int_{\overset{\quad\!\!\!\!\!\!\!\! \circ}{p_\omega^c} \times S_k} \chi_E(\omega.(x,\theta)) d\lambda(x ,\theta) \right) d\mu^t(\omega),
\end{eqnarray}
where $p^c_\omega$ is the prototile type of the tile in $\omega$ surrounding the origin.\\
\\
As $V^k_i$ and $V^l_j$ are disjoint if $i \neq j$ or $k \neq l$ and because the border of $V^k_i$ is a set of measure zero, we have : 
$$\hspace{-0.5cm} 1=\mu(\Omega) = \sum_{i=1}^{108} \sum_{k=0}^2 \mu(F^k_i) = \sum_{i=1}^{108} \mu(U_i) = \sum_{i=1}^{108} \mu^t \big( \Xi(p_i^c) \big) = \sum_{i=1}^{108} \mu^t \big ( \chi_{p_i^c} \big ) = \mu^t (\Xi)$$

		\end{enumerate}
\begin{flushright}
$\square$
\end{flushright}

\noindent
Hence, we have obtained the following result :

\begin{prop}
For pinwheel tilings, the gap-labelling (or patch frequencies) is given by :
$$\tau^\mu_* \Big( K_0 \big(C(\Omega) \rtimes \RR^2 \rtimes S^1 \big) \Big) = \frac{1}{264} \ZZ \left[ \frac{1}{5} \right].$$
\end{prop}

\noindent
The second point is that we can recover from the explicit formula (\refg{pp}) for $\mu^t$ that, in fact, the continuous hull of the pinwheel tiling equipped with the action of the direct isometries is uniquely ergodic :

\begin{lem}
$\big( \Omega,\RR^2 \rtimes S^1 \big)$ is uniquely ergodic.
\end{lem}

\noindent
\textbf{Proof :}
		\begin{enumerate}
			\item[] This is readily obtained because $\Xi$ is a Cantor set and hence every Borel measure on $\Xi$ is completely determined by its value on each clopen partition of $\Xi$.\\
If $\mu$ and $\nu$ are two invariant Borel probability measures on $\Omega$ then they induce two invariant probability measures $\mu^t$ and $\nu^t$ on the transversal $\Xi$.\\
By the equality (\refg{pp}), these measures must agree on the clopen sets defined by the successive collared supertiles which form a base of clopen neighborhoods for $\Xi$ so $\mu^t = \nu^t$.\\
Finally, by formula (\refg{mu}), we can see that $\mu=\nu$.

		\end{enumerate}
\begin{flushright}
$\square$
\end{flushright}

\newpage

\section{Conclusion}

Using methods developed in \citeg{BelSav}, we thus proved the gap-labelling for the pinwheel tiling (the $(1,2)$-pinwheel tiling).\\
We think that the methods developed in the present paper are in fact more general.\\

\noindent
We think that, using the construction of Bellissard and Savinien \citeg{BelSav}, we may prove that the top \v{C}ech cohomology group of more general tilings is in fact the integer group of coinvariants of the top $\Delta$-transversal.\\
This result coupled with the diagram of section \textbf{4.1} can then prove the known fact (see \citeg{BelBenGam}) that the image under the Ruelle-Sullivan map of the copy of the top integer  \v{C}ech cohomology of the hull in the longitudinal cohomology of $\Omega$ is in fact $\mu^t(C(\Xi,\ZZ))$.\\

\noindent
We point out another possible generalization of our work.\\
We guess that, for all $m$ and $n$ in $\NN^*$, the number of $(m,n)$-pinwheel tilings fixed by a finite rotation always remains finite allowing us to apply the results of \citeg{Hai}.\\
The result of the final section can still be applied to such tilings, obtaining that the gap-labelling of $(m,n)$-pinwheel tilings is given by the $\ZZ$-module of "patch frequencies" $c.\ZZ \left[\frac{1}{m^2 + n^2} \right]$, where $c$ is a scalar normalizing the Perron eigenvector.\\

\noindent
Moreover, we think that for any substitution (self-similar, repetitive, of finite type and non periodic) tiling of $\RR^d$ with expansion scalar $\lambda$, the explicit computation can be done in the same way, showing that the gap-labelling for such a tiling is given by $c.\ZZ \left[ \frac{1}{\lambda^d} \right]$ if $\lambda$ is an integer and by 
$$c.\left\{  \frac{1}{\lambda^{d.n}} \sum \alpha_i k_i ; n \in \NN, k_i \in \ZZ \right\} $$
if $\lambda$ is not an integer, where $\alpha$ is the unique normalized (Perron-Frobenius) eigenvector of the substitution matrix of the tiling for collared prototiles.\\
This result would generalize a result proved by Kellendonk in \citeg{Kelcoinv}.\\
\\
This result would also allow us to retrieve the fact that the tiling space of any substitution tiling is uniquely ergodic, which was already known (see \citeg{Dirk}, \citeg{LeeMooSol}, \citeg{Rad1} and \citeg{Solo}).

\newpage

\newpage

\bibliographystyle{alpha}
\bibliography{article-pimsner-cohom.bib}

\text{ }\\
\\
Univ. Blaise Pascal, Clermont-Ferrand, FRANCE\\
\textit{E-mail address :}  haija.moustafa@math.univ-bpclermont.fr\\
\textit{URL :} http://math.univ-bpclermont.fr/~moustafa

\end{document}